\documentclass[aap,MSNbibl,nameyear,dvips]{arximspdf}
\usepackage{graphicx}
\doi{10.1214/11-AAP793}
\volume{22}
\issue{4}
\pubyear{2012}
\firstpage{1328}
\lastpage{1361}

\makeatletter
\newtheorem{theorem}{Theorem}
\newtheorem{lemma}{Lemma}
\newproclaim{rmk}{Remark}

\newcommand{\bbe}{\mathbb{E}}
\newcommand{\bbp}{\mathbb{P}}
\newcommand{\cT}{\mathcal{T}}
\newcommand{\cS}{\mathcal{S}}
\newcommand{\cR}{\mathcal{R}}
\newcommand{\cN}{\mathcal{N}}
\newcommand{\cM}{\mathcal{M}}
\newcommand{\cG}{\mathcal{G}}
\newcommand{\cF}{\mathcal{F}}
\newcommand{\cE}{\mathcal{E}}
\newcommand{\cD}{\mathcal{D}}
\newcommand{\cC}{\mathcal{C}}
\newcommand{\cB}{\mathcal{B}}
\newcommand{\cA}{\mathcal{A}}
\renewcommand{\P}{\mathbb{P}}
\newcommand{\E}{\mathbb{E}}
\newcommand{\eps}{\varepsilon}
\newcommand{\var}{\operatorname{Var}}
\newcommand{\cov}{\operatorname{Cov}}

\makeatother

\begin{document}
\begin{frontmatter}

\title{Effect of scale on long-range random graphs and chromosomal inversions}
\runtitle{Effect of scale}

\begin{aug}
\author[A]{\fnms{Nathana\"el} \snm{Berestycki}\ead[label=e1]{N.Berestycki@statslab.cam.ac.uk}}
\and
\author[A]{\fnms{Richard} \snm{Pymar}\corref{}\ead[label=e2]{R.Pymar@statslab.cam.ac.uk}}
\runauthor{N. Berestycki and R. Pymar}
\affiliation{University of Cambridge}
\address[A]{Statistical Laboratory\\
University of Cambridge\\
Wilberforce Road\\
Cambridge, CB3 0WB\\
United Kingdom\\
\printead{e1}\\
\hphantom{E-mail: }\printead*{e2}} %adresu isvedimo komanda gale!
\end{aug}

% HISTORY:
\received{\smonth{3} \syear{2011}}
\revised{\smonth{6} \syear{2011}}

% ABSTRACT
%
\begin{abstract}
We consider bond percolation on $n$ vertices on a circle where edges
are permitted between vertices whose spacing is at most some number
$L=L(n)$. We show that the resulting random graph gets a~giant
component when $L\gg(\log n)^2$ (when the mean degree exceeds~1) but
not when $L\ll\log n$. The proof uses comparisons to branching random
walks. We also consider a related process of random transpositions of
$n$ particles on a circle, where transpositions only occur again if the
spacing is at most $L$. Then the process exhibits the mean-field
behavior described by Berestycki and Durrett if and only if $L(n)$
tends to infinity, no matter how slowly. Thus there are regimes where
the random graph has no giant component but the random walk
nevertheless has a phase transition. We discuss possible relevance of
these results for a dataset coming from D. repleta and D. melanogaster
and for the typical length of chromosomal inversions.
\end{abstract}

% KEYWORDS
%
\begin{keyword}[class=AMS]
\kwd[Primary ]{05C80}
\kwd{60K35}
\kwd[; secondary ]{92D15}.
\end{keyword}
\begin{keyword}
\kwd{Random transposition}
\kwd{random graphs}
\kwd{phase transition}
\kwd{coagulation-fragmentation}
\kwd{giant component}
\kwd{percolation}
\kwd{branching random walk}
\kwd{genome rearrangement}.
\end{keyword}

\end{frontmatter}

%s1 #&#
\section{Introduction and results}\label{intro}

%s1.1 #&#
\subsection{Random graphs results}\label{sec1.1}

Let $n\ge1$ and let $L=L(n) \ge1$. Define vertex set $V= \{1, \ldots,
n\}$ and edge set $\mathcal{R}_L = \{(i,j) \in V^2, \|i-j\|\le L\}$,
where $\| i - j\|$ denotes the cyclical distance between $i$ and $j$,
that is, $\|u\|= \min(|u|,\allowbreak n-|u|)$ for $u \in V$. In this paper we
consider bond percolation on $V$ where each edge in $\cR_L$ is open
with probability $p$. Equivalently, let $(G(t),t\ge0)$ be the random
graph process where a uniformly chosen edge of $\cR_L$ is opened in
continuous time, at rate 1. Let $\Lambda^1(t) \ge\Lambda^2(t)\ge
\cdots$ denote the ordered component sizes of $G(t)$. At a fixed time
$t$ this corresponds to the above model with $p = 1- \exp\{ - t/(nL)\}
$. When $L=1$, this is the usual bond percolation model on the cycle of
length $n$, while for $L(n)=n/2$, we find that $G(t)$ is a realization
of the much studied random graph model of Erd\H{o}s and Renyi [see
\citet{bollobasu} and \citet{durrett} for background]. Hence,
our random graph model interpolates between these two cases.

In this paper we are interested in the properties of the connected
components of $G(t)$, particularly those related to the possible
emergence of a~giant component when the average degree exceeds 1. The
main result of this paper shows that this depends on the scale $L(n)$.
To state our results, we let $c>0$ and consider $t = cn/2$, so that the
expected degree of a given vertex in~$G(t)$ converges to $c$ when $n
\to\infty$. Let $\Lambda^1(t) \ge\Lambda^2(t)\ge\cdots$ denote the
ordered component sizes of $G(t)$.
%
%th1 #&#
\begin{theorem}\label{TuniqGC}
Let $t=cn/2$, where $c>0$ is fixed as $n \to\infty$.

\begin{longlist}
\item If $c<1$, then there exists $C<\infty$ depending only
on $c$ such that $\Lambda^1(t) \le C \log n$ with high probability as
$n \to\infty$.

\item If $c>1$ and there exists $\xi>0$ such that $L(n)
\ge(\log n)^{2+\xi}$, then there is a unique giant component; more precisely,
%
%e1 #&#
\begin{equation}
\label{perc}
\frac{\Lambda^1(t)}{n} \to\theta(c)
\end{equation}
in probability as $n \to\infty$, where $\theta(c)$ is the survival
probability of a $\operatorname{Poisson}(c)$ Galton--Watson tree. Moreover, $\Lambda
^2(t)/n \to0$ in probability.

\item However, if $c>1$ and $L=o(\log n)$, then for all $a>0$,
%
%e2 #&#
\begin{equation}
\label{noperc}
\frac{\Lambda^1(t)}{n^a} \to0
\end{equation}
in probability as $n\to\infty$. In particular there are no giant components.
\end{longlist}
\end{theorem}

Statement (i) is fairly easy to prove using the standard technique of
approximating the size of a~component in the graph by the total progeny
of a~branching process. The result follows since in the case $c<1$ we
know that the total progeny of the branching process is almost surely
finite and has exponential tails.

Part (ii) is the most challenging. We start by noting that the
exploration of the component containing a given vertex $v$ may be
well-approximated by the trace of a branching random walk where the step
distribution is uniform on $\{- L , \ldots, L\}$. This approximation is
valid so long as the local density of the part of the component already
explored stays small. Thus, showing the existence of a giant component
requires a balancing act; we need to ensure that the local density of
what we explore stays small enough to ignore self-intersections, but
large enough for global connections to occur. Careful estimates on
survival probabilities of killed branching random walks are used to
achieve this.

Part (iii) is the easiest to prove, and requires showing the existence
of many ``blocking'' intervals of size $L$ which consist just of
vertices with degree~0. When there are many such intervals, no giant
component can exist.

%s1.2 #&#
\subsection{Long-range random transpositions}

Theorem \ref{TuniqGC} was originally motivated by the study of a
question concerning long-range transpositions, which may itself be
rephrased as a question in computational biology. We now discuss the
question on long-range random transpositions and delay the applications
to comparative genomics until Section \ref{Sgenomics}.

Recall the definitions of $V$ and $\cR_L$ in Section \ref{sec1.1}. Consider a
random process $(\sigma_t,t\ge0)$ with values in the symmetric group
$\cS_n$, which evolves as follows. Initially, $\sigma_0= e$ is the
identity permutation. Let $(i_1, j_1), (i_2, j_2), \ldots$ be an i.i.d.
infinite sequence of pairs of elements of $V$, where each pair is
uniformly distributed on $\cR_L$. Then we put
%
%e3 #&#
\begin{equation}
\sigma_t = \tau_{N_t} \circ\cdots\circ\tau_1,
\end{equation}
where for each $k\ge1$ we let $\tau_k$ denote the transposition
$(i_k,j_k)$, \mbox{$(N_t,t\ge0)$} is an independent Poisson process with rate
1 and $\circ$ the composition of two permutations. That is, informally,
if we view the permutation $\sigma_t$ as describing the positions on
the circle of $n$ particles labeled by $V$ [with $\sigma_t(i)$ denoting
the position of particle $i \in V$], then in
continuous time at rate~1, a~pair of positions $(i,j)$ is sampled
uniformly at random from $\mathcal{R}_L$ and the two particles at
positions $i$ and $j$ are swapped. Thus the case where $L(n)\ge n/2$
corresponds to the well-known random transposition process (i.e., the
composition of uniform random transpositions), whereas the case where
$L(n) =1$ corresponds to the case of random adjacent transpositions on
the circle.

Our interest consists of describing the time-evolution of $\delta
(\sigma
_t)$, where for all $\sigma\in S_n$ we set $\delta(\sigma) = n -
|\sigma|$ and $|\sigma|$ to be the number of cycles of $\sigma$. By a
well-known result of Cayley, this is the length of a minimal
decomposition of $\sigma$ into a product of \textit{any} transpositions
(i.e., whose range is not necessarily restricted to~$\mathcal{R}_L$).
The reason for this choice will become apparent in subsequent sections
and is motivated by the applications to comparative genomics.

For $c>0$, define a function
%
%e4 #&#
\begin{equation}\label{Du}
u(c)=1-\sum_{k=1}^{\infty}\frac{1}{c}\frac{k^{k-2}}{k!}(ce^{-c})^k.
\end{equation}
It is known that $u(c) = c/2$ for $c\le1$ but $u(c)<c/2$ for $c>1$
[see, e.g., \citet{bollobasu}, Theorem 5.12]. The function $u$ is
continuously differentiable but has no second derivative at $c=1$. We
shall prove the following results.

%th2 #&#
\begin{theorem}\label{Tdist}
Assume $L(n) \to\infty$ as $n \to\infty$. Then we have the following
convergence in probability as $n\to\infty$: for all $c>0$,
%
%e5 #&#
\begin{equation}\label{Eu}
\frac1n \delta( \sigma_{cn/2}) \to u(c).
\end{equation}
\end{theorem}

%Roughly speaking, this result says that up to time $n/2$ the number of
%cycles decreases linearly at rate 1 and thus only coagulations occur
%(no fragmentations). After this time however the rate that cycles
%decrease is less than 1 which suggests that fragmentations do occur
%and hints at the existence of non-microscopic cycles.

In this result the distance between the two points being transposed at
every transposition is uniform within $\{1, \ldots, L(n)\}$. We will
prove in Theorem~\ref{Tdistgeneral} given in Section \ref{Sproofofbio} a more general version of this result, where this length is
allowed to be some arbitrary distribution subject to the condition that
there are no ``atoms in the limit,'' which is the equivalent of
requiring here $L(n) \to\infty$.

By contrast, the microscopic regime (where $L$ is assumed to be
constant or to have a limit) shows a remarkably different behavior.
%
%th3 #&#
\begin{theorem}\label{Tmicro}
%Suppose $\limsup_{n\to\infty} L(n)<\infty$. Then with probability
%$1-o(1)$, for all $c>0$, $$\lim_{n\to\infty}\frac{1}{n}\delta(
Assume $\lim_{n\to\infty}L(n)$ exists. Then we have convergence in
probability; for all $c>0$,
\[
\frac{1}{n}\delta(\sigma_{cn/2})\to v(c)
\]
as $n\to\infty$, for some $C^2$ function $v(c)$ which satisfies
$0<v(c)<c/2$ for all $c>0$.
\end{theorem}

%In particular Theorem \ref{Tmicro} says that if $\lim_{n\to\infty}
%L(n)$ exists then fragmentations occur at all times. This is
%relatively simple to prove (see Lemma \ref{Lboundedfrags}) by showing
%that for any $c>0$ by time $t=cn/2$ there is a positive probability
%that any allowed 2-cycle has been created and then fragmented.

As we will describe in greater detail later on, there is a connection
between the long-range random transposition process and the random
graph process of Theorem \ref{TuniqGC}. Roughly speaking, when $L(n)$
is bounded, we expect $v(c)<c/2$ because each new edge has a positive
probability of having its two endpoints in the same connected
component. Alternatively, the branching random walk which is used to
explore the connected component of a vertex has a positive probability
of making a self-intersection at every new step.

The mean-field case where $L(n) = n/2$ recovers Theorem 4 of
\citet{beres}. Theorem \ref{Tdist} above relies on a coupling with the
random graph $G(t)$ of Theorem \ref{TuniqGC}; this coupling is similar
to the coupling with the Erd\H{o}s--Renyi random graph introduced in
\citet{beres}. In that paper, the emergence of the giant component in
the Erd\H{o}s--Renyi random graph was a crucial aspect of the proofs.
As a result, one might suspect that the phase transition of $\delta
(\sigma_t)$ is a direct consequence of the emergence of a giant
component in the random graph. However, one particularly surprising
feature of Theorem \ref{Tdist} above is the fact that the limiting
behavior described by (\ref{Eu}) holds for all $L(n) \to\infty$, no
matter how slowly. This includes in particular the cases where $L(n) = o (
\log n)$ and the random graph~$G(t)$ does \textit{not} have a giant
component. Hence, for choices of $L(n)$ such that $L(n) \to\infty$ but
$L(n) = o(n)$, the quantity $\delta(\sigma_t)$ has a phase transition
at time $n/2$, even though the random graph $G(t)$ does not get a giant
component at this time.

%s1.3 #&#
\subsection{Relation to other work, and open problems}

\textit{Long-range percolation}. A~similar model has been studied by
\citet{penrose}. There the model considered is on the \textit{infinite} square
grid $\mathbb{Z}^d$, rather than the finite (one-dimensional) torus
which we consider here. In the \textit{infinite} case, $d=1$ is trivial
since percolation (occurrence of an infinite cluster) only occurs if
$p=1$ for obvious reasons. Penrose studied the case $d\ge2$ and showed
that if $c$ is the expected degree of the origin, and $L$ the maximum
distance between the two ends of a bond, where the parameter $L\to
\infty
$ and $c$ is fixed, then the percolation probability approaches $\theta
(c)$ [where $\theta(c)$ is the same as in (\ref{perc}), i.e., the
survival probability for a Galton--Watson process with Poisson($c$)
offspring distribution]. As is the case here, his arguments use a
natural comparison with branching random walks.

It is interesting that, while the infinite case is essentially trivial
when $d=1$, the \textit{finite}-$n$ case is considerably more intricate
than the infinite case, as witnessed by the different behaviors in
(\ref{perc}) and (\ref{noperc}) depending on how fast $L(n) \to
\infty
$. Regarding the finite-$n$ situation, it is an interesting open
question to see whether there are giant components if $t=cn/2$ and
$c>1$ with $\log n \le L(n) \le(\log n)^2$. Another interesting
problem concerns the size of the largest components when there is no
giant component, in particular, if $L =o(\log n)$. Indeed, our
proof makes it clear that when $L(n) \to\infty$, even if the largest
component is not macroscopic, there is a positive proportion of
vertices in components of \textit{mesoscopic} size.
We anticipate that as $c>1$ is fixed and $L$ increases, the size of the
largest component, normalized by $n$, jumps from 0 to $\theta(c)$ as
$L(n)$ passes through a critical threshold between $\log n$ and $( \log n)^2$.
As pointed out by a referee, this is suggested by a work of
\citet{aizen} on long-range bond percolation on $\mathbb{Z}$ where the
connection probability between vertices at distance $x>0$ decays like
$1/ x^2$. Their main result (Proposition 1.1) shows that such
discontinuities occur in this case.

\textit{Epidemic models.} The question of giant components in random
graph models can, as usual, be rephrased in terms of epidemic
processes. More precisely, fix a vertex $v$ and a number $p \in
(0,1)$. Consider an SIR epidemic model that begins with all vertices
susceptible but vertex $v$ infected. Once a vertex is infected, it
transmits the infection to each of its neighbors in the base graph
$(V,E)$ at rate $\lambda>0$ and dies or is removed at rate 1. Then the
total size of the epidemic is equal to the size of the component
containing~$v$ in the random graph with edge-probability $p = \lambda/
(1+ \lambda)$. As pointed out by an anonymous referee, \citet{bramson}
consider the related SIS model (or contact process) on $\mathbb{Z}^d$
where, as here, long-range connections are possible. Similar
techniques are employed as in this article to calculate the critical
rate of infection and the probability of percolation. Letting
infections occur at rate $\lambda/ \operatorname{Vol} B(L)$ where $B(L)$ is a
ball or radius~$L$ in $\mathbb{Z}^d$, they show that the critical
infection rate $\lambda_c$ converges to~1 in all dimensions as
$L\to\infty$. They also identify the rate of convergence, which turns
out to depend on the dimension in an interesting way.

\textit{Higher-dimensional analogs of Theorem} \ref{TuniqGC}. Our
proofs do not cover the higher-dimensional cases but it would not be
very difficult to adapt them. In particular, the analogue of (\ref
{perc}) would hold if $d \ge2$ no matter how slowly $L(n) \to\infty$.
In other words, only for the one-dimensional case is it important to
have some quantitative estimates on $L(n)$. Intuitively this is
because, in one dimension, one is forced to go through potentially bad
regions whereas this problem does not arise in higher dimensions.

Regarding site percolation, we point out that recently \citet{bollobas}
have described an interesting behavior for a site percolation model on
the torus in dimensions $d\ge2$ where two vertices are joined if they
agree in one coordinate and differ by at most $L$ in the other. For
$d=2$ they show that the critical percolation probability, $p_c(L)$,
satisfies $\lim_{L\to\infty}L p_c(L)=\log(3/2)$. This is surprising as
the expected degree of a~given vertex at the phase transition is then
strictly greater than 1. There again, approximation by branching random
walks plays an important role in the proof.

\textit{Slowdown transitions for random walks}. In the mean-field case
$L(n) = n/2$ of uniformly chosen random transpositions, the quantity
$\delta(\sigma_t)$ may be interpreted as the graph-theoretic distance
between the starting position of the random walk ($\sigma_0=$ the
identity element) and the current position of the walk. Theorem \ref
{Tdist} in this case [which, as already mentioned, is Theorem 4 of
\citet{beres}], may thus be interpreted as a \textit{slowdown} transition
of the evolution of the random walk; at time $n/2$, the acceleration
[second derivative of $\delta(\sigma_t)$] drops from 0 to $-\infty$. By
contrast, \citet{limit} studied the evolution of the graph-theoretic
distance in the case of random adjacent transpositions. This
essentially corresponds to the case $L=1$, with the difference that the
transposition $(1\ n)$ is not allowed. They found that no sudden
transition occurs in the deceleration of the random walk. It would be
extremely interesting to study the evolution of the graph-theoretic
distance of the random walk when $L=L(n)$ is a given function that may
or may not tend to infinity as $n \to\infty$. Unfortunately, this
problem seems untractable at the moment as it is far from obvious how
to compute (or estimate) the graph distance between two given
permutations. [We note that even in the case $L=1$ where the
transposition $(1\ n)$ is allowed, this question is partly open; see
Conjecture 3 in \citet{limit}.] Nevertheless it is tempting to take
Theorems \ref{Tdist} and \ref{Tmicro} as an indication that a
slowdown transition for the random walk occurs if and only if $L(n) \to
\infty$, with the phase transition always occurring at time $n/2$.

\textit{Organization of the paper.} In Section \ref{Sgenomics} we show
how Theorems~\ref{Tdist} and~\ref{Tmicro} relate to a biological
problem and in particular discuss the possible relevance of these
results for a dataset coming from two Drosophila species.
In Section~\ref{Scomponents} we state and prove results on the
evolution of the clusters in a~random graph which evolves in a more
general way to $G(t)$. In Section \ref{Sproofofmain} we give a~proof
of Theorem \ref{TuniqGC}.\vadjust{\goodbreak} Section~\ref{Sproofofbio} contains a
proof of a result stronger than Theorem \ref{Tdist} using the more
general random graph process defined in Section \ref{Scomponents}.
Finally, in Section \ref{Sproofofmicro} we present the proof of
Theorem \ref{Tmicro}.

%s2 #&#
\section{Applications in comparative genomics}
\label{Sgenomics}
%s2.1 #&#
\subsection{Statement of problem and history}
Part of the motivation for this paper comes from a biological
background, more specifically, in answering a~question about the
evolution of the gene order of chromosomes. We begin with an example.
In 2001 Ranz, Casals, and Ruiz located 79 genes on chromosome 2 of
Drosophila repleta and on chromosome arm 3R of Drosophila melanogaster.
While the genetic material is overall essentially identical, the order
of the genes is quite different. If we number the genes according to
their order in D. repleta then their order in D.~melanogaster is
given in Table~\ref{table}.

%t1 #&#
\begin{table}
\tablewidth=310pt
\caption{Order of the genes in D. repleta compared to their
order in D. melanogaster}
\label{table}
\begin{tabular*}{\tablewidth}{@{\extracolsep{\fill}}r r r r r r r r r r@{}}
\hline
36&37&17&40&16&15&14&63&10&9\\
55&28&13&51&22&79&39&70&66&5\\
6&7&35&64&33&32&60&61&18&65\\
62&12&1&11&23&20&4&52&68&29\\
48&3&21&53&8&43&72&58&57&56\\
19&49&34&59&30&77&31&67&44&2\\
27&38&50&26&25&76&69&41&24&75\\
71&78&73&47&54&45&74&42&46\\
\hline
\end{tabular*}
\end{table}

Since the divergence of the two species, this chromosome region has
been subjected to many reversals or chromosomal inversions, which are
moves that reverse the order of whole gene segments. Because they
involve many base pairs at a time rather than the more common
substitutions, insertions and deletions, these mutations are called
large-scale. They are usually called inversions in the biology
literature, but we stick with the word reversal as ``inversions'' is
often used among combinatorists with a different meaning [see, e.g.,
\citet{DiaconisGraham}]. One question of interest in the field of
computational biology is the following: \textit{How many such reversals
have occurred}?

\citet{hp} have devised a widely used algorithm which computes the
\textit{parsimony distance}, the minimal number of reversals that are needed
to transform one chromosome into the other (this will be denoted here
by $d_\infty$). By definition, the number of reversals that did occur
is at least $d_\infty$. \citet{beres} complemented this by rigorously
analyzing the limiting behavior of the discrepancy between the true
distance and the parsimony distance [described by the function $u(c)$
in Theorem \ref{Tdist}], under the mean-field assumption that all
reversals are equally likely.

However, that assumption does not seem to be entirely justified and it
might be more accurate to restrict the length of the segment being
reversed. According to \citet{chromosomes}, ``To seek a biological
explanation of the nonuniformity we note that the gene-to-gene pairing
of homologous chromosomes implies that if one chromosome of the pair
contains an inversion that the other does not, a loop will form in the
region in which the gene order is inverted$\ldots.$ If a recombination
occurs in the inverted region then the recombined chromosomes will
contain two copies of some regions and zero of others, which can have
unpleasant consequences. A simple way to take this into account is$\ldots$
[to] restrict our attention to the $L$-reversal model.'' The
reasoning here is that as the length of the segment reversed increases,
the probability of recombination increases. Here, the $L$-reversal
model is to allow only reversals that switch segments of, at most,
length $L$ and all such reversals have equal probability. A further
argument can be seen in \citet{DNA} who argues that not all inversions
occur at the same rate; when a~large amount of DNA is absent from a
chromosome, the offspring is typically not viable, so longer inversions
will occur at a lower rate.

%s2.2 #&#
\subsection{Estimating the number of chromosomal inversions}

To estimate the number of chromosomal inversions (or reversals) in the
long-range spatial model, one natural idea is to use the parsimony
approach; that is, compute the $d_L$-distance (minimal number $d_L$ of
$L$-reversals needed to transform one genome into the other) and then
prove a limit theorem for the evolution of $d_L(t)$ under random
$L$-reversals. However, this appears completely out of reach at this
stage; the crucial problem is that we do not know of any algorithm to
compute the $L$-reversal distance. [Even in the case $L=1$, if
particles are lying on a circle, this is a delicate problem; see
Conjecture 3 in \citet{limit}.] Thus, even if a limit theorem could be
proved, we would not know how to apply it to two given genomes.

In order to tackle this difficulty, we propose here the following
alternative approach. We keep looking at the $d_\infty$-distance
(minimal number of reversals needed to transform one chromosome into
the other, no matter their length) but now we think of $d_\infty$ only
as an easily computed statistic on which we can make some inference,
even though not all reversals were equally likely. More precisely, we
are able to describe the evolution of the quantity $d_\infty(t)$ under
the application of random $L$-reversals, and use that result to
estimate~$t$ from the data $d_\infty(t)$.

We first state the result in this context, and illustrate our idea with
a~numerical example in Section \ref{SSnum}. The distance $d_\infty$ is
defined in terms of an object known as the \textit{breakpoint graph}.
For definitions of these notions we refer the interested reader to
Chapter 9 of \citet{DNA}. For signed permutations $\sigma,
\sigma'$ we let $\hat \delta(\sigma, \sigma') = n+ 1 - c$, where $c$ is
the number of components of the\vadjust{\goodbreak} breakpoint graph. In
general [see \citet{DNA}, Theorem~9.1], $d_\infty(\sigma, \sigma')
\ge\hat\delta(\sigma, \sigma ')$. The quantity $\hat\delta(\sigma,
\sigma')$ ignores obstacles known as ``hurdles'' and ``fortresses of
hurdles.'' All these are thought to be negligible in biologically
relevant cases, so we will use $\hat \delta(\sigma, \sigma')$ as a
proxy for $d_\infty(\sigma, \sigma')$. Let~$\sigma_t$ be the signed
permutation obtained by composing $\operatorname{Poisson}(t)$
independent $L$-reversals. We slightly abuse notation and write
$\hat\delta(\sigma_t)$ for $\hat\delta(\sigma_0, \sigma_t)$.

%th4 #&#
\begin{theorem}\label{Tdist2}
Assume that $L(n) \to\infty$. Then
\[
\frac1n \hat\delta( \sigma_{cn/2}) \to u(c).
\]
\end{theorem}

However, when $L(n)$ stays bounded, we get a behavior similar to
Theorem~\ref{Tmicro}.
%
%th5 #&#
\begin{theorem}\label{Tmicro2}
Assume $\lim_{n\to\infty}L(n)$ exists. Then we have convergence in
probability; for all $c>0$,
\[
\frac{1}{n} \hat\delta(\sigma_{cn/2})\to w(c)
\]
as $n\to\infty$, for some $C^2$ function $w(c)$ which satisfies
$0<w(c)<c/2$ for all $c>0$.
\end{theorem}

The proofs for these two results are \textit{verbatim} identical to those
of Theorems \ref{Tdist} and \ref{Tmicro}. The choice of stating our
results for transpositions is merely one of convenience, as
transpositions are easier to describe and more familiar to many mathematicians.

%s2.3 #&#
\subsection{Numerical application to Drosiphila set}\label{SSnum}

We now illustrate on the dataset from Table \ref{table} the possible
relevance of Theorems \ref{Tdist2} and \ref{Tmicro2}. We first
compute the parsimony distance in this case. Here there are $n=79$
genes, and even though the orientation of each gene is not written, it
is not difficult to find an assignment of orientations which minimizes
the parsimony distance $d_\infty(\sigma)$. We find that the parsimony
distance is $d_\infty(\sigma) =54$.

First assume that all reversals are equally likely, or that $L$ is
large enough that the behavior described in Theorem \ref{Tdist2}
holds, and let us estimate the actual number of reversals that were
performed. We are thus looking for~$t$ such that $d_\infty(\sigma_t) =
54$ when $n=79$. Changing variables $t=cn/2$, we are looking for $c>0$
such that $u(c) = 54/79 \approx0.68$. Thus, inverting $u$ we find
$c\approx1.6$ and hence, we may estimate the number of reversals to be
around $t=63$. Note that the discrepancy with parsimony ($d_\infty=
54$) is already significant.

This estimate keeps increasing as $L$ decreases and the behavior of
Theorem \ref{Tmicro} starts kicking in. For instance, with $L=4$ (so
that $L/n \approx5\%$), simulations give $c \approx2.4$ or $t \approx
95 $ reversals, or $175\%$ of the initial parsimony estimate!

%f1 #&#
\begin{figure}[b]
\vspace*{-3pt}
\begin{tabular}{@{}ccc@{}}

\includegraphics{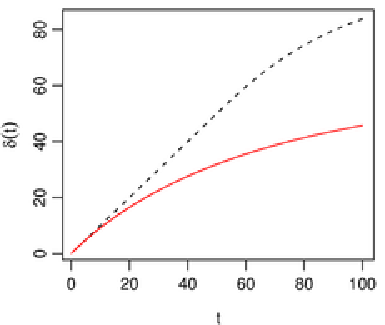}
 & \includegraphics{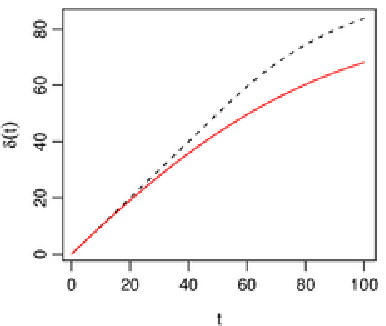} & \includegraphics{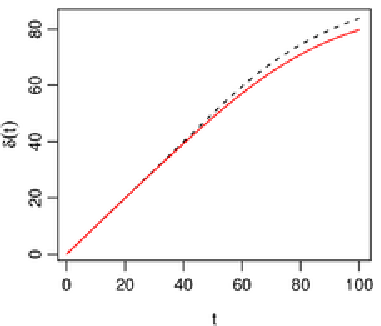}\\
(a) & (b) & (c) \\[5pt]

\includegraphics{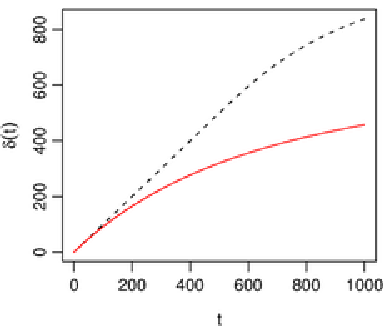}
 & \includegraphics{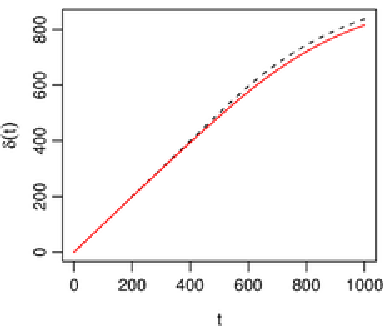} & \includegraphics{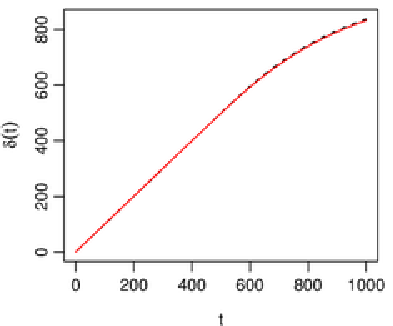}\\
(d) & (e) & (f)
\end{tabular}
\caption{Here \textup{(a)} $n=100, L=1$; \textup{(b)} $n=100, L=5$;
\textup{(c)} $n=100, L=50$; \textup{(d)}~$n=1000,\allowbreak L=1$;
\textup{(e)} $n=1000, L=50$; \textup{(f)} $n=1000, L=500$.}
\label{fig1}
\end{figure}

Ideally, we would want to use estimates in the biology
literature on the typical range of reversals, in combination with the
results of this paper, to produce a refined estimate. \citet{haussler}
estimated the median length of a reversal in human/mouse genomes to be
about $1$~kb, corresponding very roughly speaking to $L$ being a few
units, say $1\le L \le4$. (However, they find a distribution for the
reversal lengths which is bimodal and hence, quite different from the
one we have chosen for simplicity in this paper.)
Other estimates we found in several biology papers differed by several
orders of magnitude, so that there does not appear to be a consensus on
this question. Instead, we performed some reverse engineering, and
compared our method with other existing methods. \citet{YDN} used a
Bayesian approach in a model comparable to ours. The mode of the
posterior distribution was at $t\approx87$, with the parsimony
estimate lying outside the 95\% confidence interval [see \citet{DNA},
Section 9.2.2, for further details]. This suggests that $L$ is slightly
more than 4, broadly speaking consistent with the estimate of
\citet{haussler}.\vspace*{-3pt}

%s2.4 #&#
\subsection{Simulations for transpositions}
We complement the above example with plots (see Figure \ref{fig1}) to
show how $\hat\delta (t)=\hat\delta(\sigma_t)$ evolves with $t$ for
finite~$n$ by straightforward MCMC, averaging over 1000 simulations in
each case. The dotted line shows $u(c)$ and the solid line shows the
average over the simulations. We observe that as $L$ increases, $u(c)$
provides a better estimate to the parsimony.\vadjust{\goodbreak}

%s3 #&#
\section{Evolution of the components of the random graph}
\label{Scomponents}
We begin by proving a few results relating to the components of a
random graph which evolves in a more general way than previously
defined. For each $n \ge2$, fix a probability distribution $(p^n_{\ell
})_{1\le\ell\le\lfloor n/2\rfloor}$. We will omit the superscript
$n$ in all calculations below in order to not overload the notation.
For the rest of this section we redefine $(G(t),t\ge0)$ to be the
random graph process where at rate 1 we choose a random variable $D$
according to the distribution $(p_{\ell})$, and open a uniformly chosen
edge from those of graph distance $D$. We define
%
%e6 #&#
\begin{equation}
\label{epsn}
\eps_n: = \max_{1\le\ell\le\lfloor n/2\rfloor} p_\ell.
\end{equation}
We begin by analyzing how the components in the random graph $G(t)$
evolve over time.
%
%le1 #&#
\begin{lemma}\label{Lbreadth}
Let $C(t)$ be the connected component of $G(t)$ containing some fixed
vertex $v \in V$, and let $t = cn/2$ for some $c>0$. Assume that $\eps
_n\to0$. We have that $C(t) \preceq Z$ (where $\preceq$ stands for
stochastic domination) and
\[
|C(t)|\to Z\qquad\mbox{as }n\to\infty
\]
in distribution. Here $Z$ is the total progeny of a Galton--Watson
branching process in which each individual has a $\operatorname{Poisson}(c)$ number of
offspring.
\end{lemma}
%
%re1 #&#
\begin{rmk}
$\!\!\!$The argument below is simpler to follow in the case where~$(p_\ell)$
is the uniform distribution in $\{1, \ldots, L(n)\}$. We will need a
more precise version of this lemma later on, in Lemma
\ref{Lbrwcoupl}. It may thus be helpful for the reader to look at Lemma
\ref{Lbrwcoupl} first in order to understand the main idea of the
proof of Lemma \ref{Lbreadth}.
\end{rmk}
\begin{pf*}{Proof of Lemma \ref{Lbreadth}}
We use the \textit{breadth-first search} exploration of the component
$C(t)$. That is, we expose the vertices that form $C(t)$ by looking
iteratively at neighborhoods of increasing radius about $v$. In doing
so, the vertices of $C(t)$ are naturally ordered according to levels
$\ell= 1, 2, \ldots$ which represent the distance of any vertex from
that level to the vertex $v$. To be more precise, if $A \subset V$ let
$\cN(A)$ denote the neighborhood of $A$, that is,
\[
\cN(A) = \{y \in V\dvtx y \sim x \mbox{ in $G(t)$ for some $x \in A$}\}.
\]
Let $A^0 = \{v\}$ and then define inductively for $i \ge0$,
\[
A^{i+1} = \cN( A^i) \Bigm\backslash\bigcup_{j=0}^i A^j.
\]
The statement of the lemma will follow from the observation that when
$t = cn/2$, the sequence $(|A^0|, |A^1|, \ldots)$ converges in the
sense of finite-dimensional distributions toward $(Z^0, Z^1, \ldots)$,
the successive generation sizes of a $\operatorname{Poisson}(c)$ Galton--Watson
process. Thus, fix an integer-valued sequence $(n_0, n_1, \ldots)$
with $n_0 =1$. We wish to show that
\[
\P(|A^0| = n_0,\ldots, |A^i| = n_i) \to\P(Z^0 = n_0, \ldots, Z^i = n_i)
\]
as $n \to\infty$,
which we do by induction on $i\ge0$. The statement is trivial for $i =
0$. Now let $i\ge0$. Given $\cA_i = \{A^0=n_0, \ldots, A^i=n_i\}$, we
look at the neighbors in level $i+1$ of each vertex $v_1,\ldots
,v_{n_i}$ in level $i$, one at a time.

Let $\bar G(t)$ be the multigraph on $V$ with identical connections as
$G(t)$, but where each edge is counted with the multiplicity of the
number of times the transposition $(i,j)$ has occurred prior to time $t$.
Equivalently, for each unordered pair of vertices $(i,j) \in V^2$ at
distance $\|i-j\|=\ell\ge1$, consider an independent Poisson process
$N^{(i,j)}(t)$ of parameter $2p_\ell/n$. Then the multigraph $\bar
G(t)$ contains $N^{(i,j)}(t)$ copies of the edge $(i,j)$, while the
graph $G(t)$ contains the edge $(i,j)$ if and only if $N^{(i,j)}(t) \ge1$.

Note that if $w \in V$, then the degree $\bar d_w$ of $w$ in $\bar
G(t)$ is
\[
\bar d_w = \sum_{\ell\ge1} \operatorname{Poisson}(2 tp_\ell/n ) =_d
\operatorname{Poisson}(c).
\]
Let $\cF_i = \sigma(A^0, \ldots, A^i)$. Conditionally on $\cF_i$, order
the vertices from $A^i$ in some arbitrary order, say $v_1, \ldots, v_{n_i}$.
Observe that
\[
A^{i+1} = \bigcup_{j=1}^{n_i} \Biggl[\cN(v_j) \Bigm\backslash\Biggl(
\bigcup
_{j=0}^i A^j \cup\bigcup_{k=1}^{j-1} \cN(v_k)\Biggr)\Biggr].
\]
It follows directly that, conditionally on $\cF_i$,
%
%e7 #&#
\begin{equation}\label{Esd}
|A^{i+1}| \preceq\sum_{j=1}^{n_i} P_j,
\end{equation}
where $P_j$ are independent $\operatorname{Poisson}(c)$ random variables which are
further independent from $\cF_i$. (The stochastic domination $|C_v|
\preceq Z$ already follows from this observation.)
For $1\le j \le n_i$, let $\cF_{i,j} = \cF_i \vee\sigma( \cN(v_1)
\cup\cdots\cup\cN(v_j))$. Observe that, conditionally given $\cF
_{i,j-1}$, then
%
%e8 #&#
\begin{equation}
N_j= \Biggl| \cN(v_j) \Bigm\backslash\Biggl( \bigcup_{j=0}^i A^j \cup
\bigcup
_{k=1}^{j-1} \cN(v_k)\Biggr) \Biggr|
\end{equation}
is stochastically dominated by $P_j$ but also dominates a thinning of
$P_j$ which is a Poisson random variable with parameter $c (1- M \eps
_n)$, where $M= n_0 + \cdots+ n_i + |\bar\cN(v_1)| + \cdots+ | \bar
\cN(v_{j-1})|$, where $\bar\cN(w)$ denotes the neighborhood of $w$ in
$\bar G(t)$ (hence, neighbors are counted with multiplicity).
Furthermore, note that the random variables $(N_j, 1\le j \le n_i)$ are
conditionally independent given $\cF_i$.
Since $\E( M\eps_n | \cF_i) \to0$ by the\vadjust{\goodbreak} stochastic domination~(\ref{Esd}), it follows that
\[
\P( |A^{i+1}| = n_{i+1} | \cF_i ) \mathbf{1}_{|A^i| = n_i} \to\P
\Biggl(\sum_{j=1}^{n_i} P_j = n_{i+1}\Biggr).
\]
%
%Thus the probability of at least one multiple edge tends to 0 as $n\to
%$G(t)$ and $\bar G(t)$ agree and thus each vertex $w \in A^i$ is
%connected to a Poisson$(\hat c_n(w))$ number of vertices which are not
%in $A^0 \cup\ldots\cup A^{i}$, where $c(1-m\eps_n/2)<\hat c_n(w)<c$
%and $m$ is the number of neighbours of $w$ in the base graph $G$ which
%are already present in $A^0 \cup\ldots\cup A^{i}$. While these
%numbers are independent for different vertices in $A^i$, the sampling
%is done with replacement and thus it could be that two vertices $w, w'
%connection as a self-intersection. %A self-intersection (of type II)
%can also occur if $w \in A^i$ selects $w' \in A^i$ as one of its
%neighbours. Given $|A^0| = n_0, \ldots|A^i|= n_i$, then the expected
%number $S^{(2)}_i$ of self-intersections of type II in level $i$ is
%given by
%$$
%$$
%(To see where this formula comes from, observe that a type II
%self-intersection can occur between any pair of the $n_i$ individuals
%in level $i$ and each possible edge exists with probability $p$.)
%We let $S_i$ be the number of self-intersections in level $i$. By a
%similar argument to above, we have

%$$
%$$
%as $n \to\infty$. Thus the probability of at least one
%self-intersection tends to 0 as $n \to\infty$ by Markov's inequality.
%However, given that there are no self-intersections and no multiple
%edges, $|A^{i+1}| = \sum_{j=1}^{n_i} P^j$, where $P^j$ are independent
%Poisson random variables with parameters $\hat c_n(j)\to c$ as $n\to

%We deduce that
%$$
%= n_{i+1}),
%$$
%where $P_j$ are independent Poisson$(c)$ random variables.
This completes the induction step and finishes the proof of convergence
in distribution.
\end{pf*}

A useful consequence of this result is the following lemma.
%
%le2 #&#
\begin{lemma}\label{LNclusterRG}
Let $t =cn/2$, where $c>0$. Then as $n \to\infty$, the number, $K_t$,
of connected components of $G(t)$ satisfies
%
%e9 #&#
\begin{equation}
\E(K_t)\sim n\sum_{k=1}^\infty\frac{k^{k-2}}{ck!}(ce^{-c})^k.
\end{equation}
\end{lemma}
\begin{pf}
For $v \in V$, let $C_v$ be the component containing vertex $v$. Then
observe that the total number of components is given by
$\sum_{v \in V}\frac{1}{|C_v|}
$
and thus by exchangeability, the expected number of components is
$n\E(1/{|C_v|})$. Dividing by $n$ and applying the bounded convergence
theorem (since \mbox{$1/|C_v| \le1$}) as well as Lemma~\ref{Lbreadth}, we obtain
\[
\frac1n \E(K_t)
\to\sum_{k=1}^\infty\frac{1}{k} \bbp(Z=k)= \sum_{k=1}^\infty
\frac
{k^{k-2}}{ck!}(ce^{-c})^k,
\]
where the exact value $\P(Z=k)$ of the probability mass function of $Z$
is the well-known Borel--Tanner distribution [see, e.g., \citet{beres},
Corollary 1].
\end{pf}

We now prove that the number of components $K_t$ is concentrated around
its mean.
%
%le3 #&#
\begin{lemma}
\label{LconcX}
Let $c>0$ and let $t= cn/2$. Assume that $\eps_n\to0$. We have $K_t/
\E
(K_t) \to1$ in probability as $n \to\infty$.
\end{lemma}
\begin{pf}
We write $K_t = Y_t + W_t$ where $Y_t$ counts the components smaller
than a threshold $T=(1/\eps_n)^{1/4}$ and $W_t$ those that are greater
than this threshold. Note that $W_t \le n/T = o(n)$ and thus it
suffices to show that $Y_t$ is concentrated around its mean, that is,
$\var(Y_t) = o(n^2)$.

Note that we can always write
\[
Y_t = \sum_{v \in V} \frac1{|C_v|} \mathbf{1}_{\{|C_v| \le T\}}
\]
and thus
\[
\var(Y_t) = n \var\bigl((1/|C_v|) \mathbf{1}_{\{|C_v| \le T\}}\bigr) + \sum
_{v \neq w}\!
\cov
\biggl(\frac1{|C_v|}\mathbf{1}_{\{|C_v| \le T\}} , \frac
1{|C_w|}\mathbf{1}_{\{|C_w| \le T\}}\biggr)\!.\vadjust{\goodbreak}
\]
Since $1/|C_v| \le1$, the first term in the right-hand side is smaller
than $n$. Define $S_v=\frac1{|C_v|}\mathbf{1}_{\{|C_v| \le T\}}$,
$S_w=\frac 1{|C_w|}\mathbf{1}_{\{|C_w| \le T\}}$. To\vspace*{2pt} know
the value of~$S_v$ and~$S_w$, it suffices to explore by breadth-first
search a relatively small number of vertices in the components of $v$
and $w$. While we do so, it is unlikely that the exploration of these
components will ever intersect, hence, the random variables $S_v$ and
$S_w$ are nearly independent.

To formalize this idea, let $C_v^T$ (resp., $C_w^T$) denote the subset
of $C_v$ (resp.,~$C_w$) obtained by exploring at most $T$ individuals
using breadth-first search as above. Let $\tilde S_v$ be a copy of
$S_v$, independent from $C_w$. Then conditionally on~$C_w^T$, exploring
$C_v$ until at most $T$ vertices have been exposed using breadth-first
search, we may take $S_v = \tilde S_v$ except if~$C_v^T$ intersects
with~$C_w^T$, an event which we denote by $\cA$. (To see this, imagine
generating an independent copy~$\tilde C_v^T$, using the same number of
offsprings and positions for each individual in the breadth-first
search exploration of $C_v^T$ as in, but stop if at any point~$\tilde
C_v^T$ has an intersection with $C_w^T$.)

Thus, letting $m_v = \E(S_v) = \E(\tilde S_v)$, since $\tilde S_v$ is
independent from $C_w^T$, and since $S_v = \tilde S_v$ on $\cA
^\complement$,
\begin{eqnarray*}
\E(S_v - m_v | C_w^T) &=& \E\bigl((S_v - \tilde S_v) | C_w^T\bigr) + \E\bigl((\tilde
S_v - m_v) | C_w^T\bigr)\\
& = &\E\bigl((S_v - \tilde S_v) \mathbf{1}_{\cA} | C_w^T\bigr)+\E\bigl((S_v -
\tilde
S_v) \mathbf{1}_{\cA^\complement} | C_w^T\bigr) \\
&=&\E\bigl((S_v - \tilde S_v)\mathbf{1}_{\cA} | C_w^T\bigr) \qquad\mbox{a.s.},
\end{eqnarray*}
and thus since $0\le S_v \le1$ and $0\le\tilde S_v \le1$,
\[
\bigl|\E(S_v - m_v | C_w^T)\bigr| \le2 \P(\cA| C_w^T)\qquad \mbox{a.s.},
\]
so that
\[
\cov(S_v,S_w) \le4\P(\cA).
\]
Now observe that by Markov's inequality, $\P(\cA) \le\E(e(C_v^T,
C_w^T))$, where\break $e(A,B)$ denotes the number of edges between $A$ and
$B$. Since $|C_v^T| \le T$ and $|C_w^T| \le T$ by definition, we have
$\E(e(C_v^T, C_w^T)) \le c\eps_nT^2/2 = O(\eps_n^{1/2}) \to0$.
The lemma follows.
\end{pf}

%s4 #&#
\section{\texorpdfstring{Proof of Theorem \protect\ref{TuniqGC}}{Proof of Theorem 1}}
\label{Sproofofmain}

%s4.1 #&#
\subsection{Connection with branching random walk}

In this section we return to considering the random graph model $(G(t),
t\ge0)$ as given in the \hyperref[intro]{Introduc-} \hyperref[intro]{tion}. The proof of
(i) in Theorem \ref{TuniqGC} is easy and follows directly from the
observation that for a given vertex $v$, $|C_v|$ is stochastically
dominated by~$Z$, the total progeny of a Poisson($c$) Galton--Watson
tree (see Lemma \ref{Lbreadth}). When $c<1$ it is easy to see that
there exists $\lambda
>0$ and $C< \infty$ such that $\P(Z>k) \le C e^{-\lambda k}$. Taking
$k=b \log n$ with $b$ sufficiently large, (i) now follows from a simple
union bound.\vadjust{\goodbreak}

We turn to the proof of (ii) in Theorem \ref{TuniqGC}, which is the
most challenging technically in this paper, and assume that $c>1$. The
key to the investigation of the properties of $G(t)$ with $t=cn/2$ is
the following observation, which connects the geometry of a given
component to the range of a certain branching random walk. We start by
introducing notation and definitions. Let $T$ be a Galton--Watson tree
with a given offspring distribution and denote by $T_i$ the $i$th level
of the tree $T$. Let $(S(v), v \in T)$ denote a $T$-indexed random
walk. That is, let $(X(e))_{e \in T}$ be a collection of i.i.d. random
variables with a prescribed step distribution, and for all vertices $v
\in T$, define $S(v) = S({o})+\sum_{e \prec v} X_v$, where the sum $e
\prec v$ runs along all edges that are on the shortest path between the
root $o$ and $v$.

Let $t=cn/2$. Let $w \in V$, say $w=0$, and let $C= C_w$ be the
component containing $w$ in $G(t)$. Consider the breadth-first
exploration of $C$ introduced in Lemma~\ref{Lbreadth}. Recall that
$A^{i+1} = \cN(A^i) \setminus\bigcup_{j=0}^i A^j$. Observe that it could
be that two vertices $w, w' \in A^i$ each select a same neighbor $z$.
We refer to this type of connection as a self-intersection. We view
each $A^i$ as a subset of $\mathbb{Z}$ by identifying $V$ with
\[
\{ - \lfloor n/2 \rfloor+1, \ldots, -1, 0, 1, \ldots, \lceil n/2
\rceil
\}.
\]
The following is a warm-up for the more complicated kind of couplings
which will be needed later on.
%
%le4 #&#
\begin{lemma}
\label{LBRW} Let $c>0$ and let $t = cn/2$. For each $k\ge1$,
\[
\biggl(\sum_{v \in A^i} \delta_{v/L}, 1\le i \le k\biggr) \to
\biggl(\sum
_{v \in T_i} \delta_{S_v}, 1\le i \le k\biggr)
\]
weakly in distribution as $n \to\infty$, where $(S_v)_{v \in T}$
denotes a branching random walk started from $0$ with offspring
distribution $\operatorname{Poisson}(c)$ and step distribution uniform on
$(-1,1)$, and $\delta_x$ denotes the Dirac pointmass at $x$.
\end{lemma}
\begin{pf}
The proof of the lemma is an easy extension of Lemma \ref{Lbreadth},
since in the case where there are no self-intersections, all the
displacements of the children of vertices in any given generation form
i.i.d. uniform random variables on $ \cN_L= \{-L, \ldots, -1, 1,
\ldots
, L\}$. Details are left to the reader.
\end{pf}

In practice, the finite-dimensional distribution convergence result of
Lem\-ma \ref{LBRW} will not be strong enough as we will typically need
to explore more than a finite number of generations. The following
lemma strengthens this to show that the breadth-first exploration of a
cluster may be coupled \textit{exactly} with a slightly modified
branching random walk up until the first time the latter has a
self-intersection.
More precisely, let $T$ be a Galton--Watson tree with offspring
distribution $\operatorname{Poisson}(c)$, and let $S_v, v \in T$, be defined as above
except that if $v \in T$ with offspring $v_1, \ldots, v_k$ [let $e_i$
denote the edge $(v,v_i)$], we define the displacement variables
$X(e_1), \ldots, X(e_k)$ to be sampled \textit{with replacement}
uniformly from~$\cN_L$. The sampling is still done independently for
different vertices $v \in T$. We call this process branching random
walk with replacement for future reference. We also introduce a version
\textit{with erasure}, where if $v$ and $w$ are such that $S_v = S_w$
(what we call a self-intersection) with~$v$ discovered before $w$ in
the breadth-first search, then the entire descendance of $w$ is ignored
or killed. We call this process an erased branching random walk and
denote it by $(\tilde S_v, v\in T)$.
%
%le5 #&#
\begin{lemma}\label{Lbrwcoupl}
Let $(S_v, v\in T)$ denote a branching random walk as above and
$(\tilde S_v, v\in T)$ its corresponding erasure. Then there exists a
coupling of $(\tilde S_v)_{v \in T}$ and $(A^i,i \ge0)$ such that the
sets $A^i$ and $\{\tilde S_v,v \in T_i\}$ coincide exactly for each $i
\ge0$. In particular, let $\tau$ be the first self-intersection level;
$\tau= \inf\{n \ge1\dvtx\allowbreak \exists v \neq w \in V(T_n), S_v = S_w\}$.
Then we can couple $A^i$ and $(S_v, v \in T_i)$ for each $i < \tau$.
\end{lemma}
\begin{pf}
For the most part this is a variation on Lemma \ref{Lbreadth}, but
there are some subtleties. Assume we are exploring the connections of a
vertex $v \in\cA^i$ for some $i\ge0$. Let $A$ be the $2L-1$ potential
neighbors of $v$, and let $B \subset A$ be the set of those within $A$
which have already been exposed so far. For each of the $|A\setminus
B|$ potential new neighbors of $v$ to be added to $\cA^{i+1}$, the edge
joining it to $v$ has appeared a $\operatorname{Poisson}(t/(nL))$ number of times. Of
course, if an edge appears several times, this amounts to connecting to
the same vertex, and this is why we choose sampling \textit{with}
replacement. The action of sampling uniformly with replacement from $A$
or from $A \setminus B$ can be chosen to be identical, until the first
time that sampling from $A$ uses an element from $B$. The rest of the
details are left to the reader.
\end{pf}
%
%re2 #&#
\begin{rmk}
Note that by the classical birthday problem, $\tau$ is unlikely to
occur before at least of order $\sqrt{L}$ vertices have been added.
Thus we can couple exactly the breadth-first search exploration of
$C_v$ and a branching random walk until of order $\sqrt{L}$ vertices
have been discovered.
\end{rmk}

In fact, this will still not be strong enough and we will need to push
this exploration until of order $o(L)$ vertices have been discovered.
Of course, self-intersections can then not be ignored, but there are
not enough of them that they cause a serious problem, so the
breadth-first search exploration coincides with ``most'' of the
branching random walk.

%s4.2 #&#
\subsection{Survival of killed branching random walk}

The basic idea for the proof of (ii) in Theorem \ref{TuniqGC} is a
renormalization (sometimes also called ``block'') argument.

We show that if the component of a given vertex is larger than some
fixed number, then this component\vadjust{\goodbreak} is likely to reach distance $KL$,
where $K>0$ is a large number (which may even depend on $L$) to be
suitably chosen. This may be iterated to show that two points selected
at random from $V$ will be connected with probability approximately
$\theta(c)^2$, where $\theta(c)$ is the survival probability of $T$.
For now, we will need a few basic estimates about killed branching
random walks. In many ways, some of the results are more natural to
state when we let $L \to\infty$ rather than $n \to\infty$. Since
$L(n) \to\infty$, the two statements are identical.

Consider a branching random walk as above, started at $v \in V$, with
step distribution uniform in $\{-L, \ldots, -1, 1, \ldots, L\}$ and
some arbitrary offspring distribution with probability generating
function $\phi(s)$. By killed branching random walk (KBRW) we refer to
a branching random walk where, in addition, particles die if they
escape a given interval containing the starting point.
%
%le6 #&#
\begin{lemma}\label{Lsurvivalprob}
Let $\theta$ denote the survival probability of the branching random
walk, that is, $\rho= 1- \theta$ is the smallest root of $z = \phi(z)$.
For each $\eps>0$ we can choose $K = K(\eps,\phi)$ such that if all
particles are killed upon escaping $[v-KL, v+KL]$, then for all $L$
sufficiently large (depending solely on $\eps$ and~$\phi$) the
survival probability $\theta^K$ of KBRW satisfies $\theta^K\ge\theta
(1-\eps)$.
\end{lemma}
\begin{pf}
Let $T$ denote the Galton--Watson tree describing the descendants
of~$v$. Conditionally on survival of $T$, the subset $U$ of $T$ for
which all vertices in $U$ have infinite progeny (i.e., the set of
infinite rays) forms a Galton--Watson process with modified progeny;
the generating function satisfies
%
%e10 #&#
\begin{equation}\label{tildephi}
\tilde\phi(s)=\frac{1}{\theta}[\phi(\theta s+1-\theta)-1+\theta].
\end{equation}
Define $\cN_v:=[v-L, v+L]$. Consider a subset $W$ of $U$ obtained as
follows. Let $\xi= (u_0, u_1, \ldots, u_R)$ be a fixed ray in $U$
where $R$ is the first time the ray leaves $[v-KL,v+KL]$. Thus
$(S_{u_0}, S_{u_1}, \ldots,S_{u_R})$ is a random walk with the
underlying step distribution, killed upon exiting $[v-KL, v+KL]$. Then
$W$ restricted to $\xi$ will consist of the subsequence of $(u_{n_i})$
such that $S_{u_{n_i}} \in\cN_v$. More precisely, we take $W$ to be
the union of all such subsequences over all rays $\xi$ in $U$. The
vertices of $W$ have a natural tree structure, and we claim that $W$
dominates a branching process where the offspring progeny is
%
%e11 #&#
\begin{equation}\label{phiK}
\phi_K (s) = \tilde\phi\bigl( \eps_K + (1-\eps_K)s\bigr),
\end{equation}
where $\eps_K = 1/(K+2)$. The reason for this is as follows. Suppose $u
\in W$, so that $S_u \in\cN_v$. Then $u$ has (in $U$) a random number,
say $N$, of offsprings, where the generating function is given by
$\tilde\phi$ in (\ref{tildephi}). Since the trajectory of a~random
walk with jumps uniform in $\{-L, \ldots, -1, 1, \ldots, L\}$ forms
a~martingale and the jumps are bounded by $L$, a classical application of
the optional stopping theorem shows that any particular fixed ray
emanating from each offspring of $u$ returns to $\cN_{v}$ before
hitting $v\pm KL$ with probability at least $1-\eps_K$. Formula (\ref
{phiK}) follows easily. Now, survival probability is therefore at least
as large as the survival probability of the Galton--Watson process with
offspring distribution given by $\phi_K$. Let $\rho_K$ be the
extinction probability. Then $\rho_K = \phi_K(\rho_K)$ and $\rho_K$ is
the unique root of this equation in $(0,1)$, and moreover, $\rho_K$ is
decreasing as a function of $K$. Since $\rho_K \ge0$, call $\rho=
\lim
_{K \to\infty} \rho_K$. It is trivial to conclude by continuity of
$\tilde\phi$ that $\rho= \tilde\phi(\rho)$ and that $\rho<1$, from
which it follows that $\rho$ is the extinction probability of $\tilde
\phi$ and is thus equal to 0. Thus we may choose $K$ sufficiently large
that $\rho_K<\eps$.\looseness=-1
\end{pf}

We now consider a certain subprocess of the killed branching random
walk and show that this also survives, growing exponentially and
leaving many offsprings very near the starting point $v$. Rather than
stating a general result we will state only what we need. Fix a
function $\omega(L$) such that $\omega(L) \to\infty$ sufficiently
slowly, say $\omega(L) = \log L$, and let $f_0(L)$ be any function such
that $f_0(L)\le L/\omega(L)$. Fix an integer $d\ge1$, and explore no
more than $d$ offsprings for any individual, that is, declare dead any
additional offspring. Fix $\lambda= \lambda_{K,d}$ and also declare a
vertex $v$ dead if the most recent common ancestor $u$ of $v$ such that
$S_u \in\cN_v$ is more than $\lambda$ generations away. Refer to this
process as $\mathit{KBRW}_1$. Note that $\mathit{KBRW}_1$ is a subprocess of $\mathit{KBRW}$ and
thus of $\mathit{BRW}$. Note also that the erased $\mathit{KBRW}_1$ is a subprocess of
the erased $\mathit{KBRW}$.
%
%le7 #&#
\begin{lemma}
\label{Lkbrw1}
Assume that $\phi''(1)< \infty$ so that the offspring distribution has
finite second moments. For all $\eps>0$, there exists $K= K(\eps,
\phi
)$, $d\ge1$ and $\lambda$ such that if all particles are also killed
upon escaping $[v-KL, v+KL]$, then for all sufficiently large $L$
(depending solely on $\phi$ and $\eps$), with probability at least
$(1-\eps) \theta$, the following hold:
\begin{longlist}
\item $\mathit{KBRW}_1$ gets at least $f_0(L)$ descendants in at most
$c \log f_0(L)$ generations for some $c>0$,
\item $f_0(L)/K$ of them are in $\cN_v$.
\end{longlist}
\end{lemma}
\begin{pf}
Consider the $\mathit{KBRW}$ of Lemma \ref{Lsurvivalprob} and let $W$ be as in
the proof of that lemma. Consider $W \cap \mathit{KBRW}_1$ and note that this is
a Galton--Watson process with offspring distribution which dominates
one with a generating function given by~(\ref{phiK}), where now
$1-\eps
_K$ is the probability that a random walk with step distribution
uniform on $\{-L, \ldots, -1, 1, \ldots, L\}$ returns to $\cN_v$ before
exiting $[v-KL,v+KL]$, and that this takes less than $\lambda$ steps.
By choosing $K$ sufficiently large, $d$ sufficiently large and $\lambda
$ sufficiently large (in that order), $\eps_K$ is arbitrarily small and
thus we find that $\mathit{KBRW}_1$ survives forever with probability at least
$(1-\eps) \theta$, as in Lemma \ref{Lsurvivalprob}. Note also that $W
\cap \mathit{KBRW}_1$, being a Galton--Watson tree and having finite second
moments, grows exponentially fast by the Kesten--Stigum theorem. Thus
fewer than $c \log f_0(L)$ levels are needed to grow $W\cap \mathit{KBRW}_1$ to
size $f_0(L)$ for some $c>0$, and so at this level we will certainly
have at least $f_0(L)$ explored in $\mathit{KBRW}_1$.

Let $\cT$ be the $\mathit{KBRW}_1$ stopped when the population size exceeds
$f_0(L)$. Define the following marking procedure in $\mathit{KBRW}_1$. Mark any
node $u \in \mathit{KBRW}_1$ if the position of the branching random walk $S_u$
at this node is in the interval $\cN_v$. Let $\cM\subset\cT$ be the
set of marked nodes. Since by construction, every node $u \in\cT$ has
an ancestor at (genealogical) distance at most $\lambda$ which is a
marked node, and since the degree of any node in $\cT$ is at most
$d+1$, it follows that
%
%e12 #&#
\begin{equation}\label{Emarked}
|\cM| \ge\eta|\cT|\qquad\mbox{where } \eta= \frac1{1+d+d^2 +\cdots+
d^\lambda}.
\end{equation}
[To see (\ref{Emarked}), just notice that for every new mark, one can
add at most $1/\eta$ nodes in the tree without adding a new mark, and
proceed by induction.] For (ii) to occur, it suffices that $|\cM| \ge
f_0(L)/K$. Since by construction $|\cT| \ge f_0(L)$, choosing $\lambda
=\lambda_{K,d} \ge\lfloor(\log K) /(\log d) \rfloor-1$ shows that
(ii) occurs as soon as (i) holds. The proof of the lemma is complete.
\end{pf}

We now strengthen this last result by showing that the erased random
walk also has a large number of offsprings in $[v-KL, v+KL]$. Further,
we suppose also that there is a set $F$ of locations which, if an
individual lands on, results in that individual being removed. We call
these \textit{forbidden locations}.
%
%le8 #&#
\begin{lemma}
\label{Lsurvival-erased}
Consider an erased branching random walk, started at \mbox{$v \in V$}, with
step distribution uniform in $\{-L, \ldots, -1, 1, \ldots, L\}$ and
some arbitrary offspring distribution with probability generating
function $\phi(s)$ with $\phi''(1)<\infty$. Suppose also that there is
a set $F$ of forbidden locations, with $|F|\le L/\omega(L)$ and
$\omega
(L)\to\infty$. Let $\theta$ denote the survival probability of the
branching random walk.
For all $\eps>0$ we can choose $K = K(\eps,\phi)$ such that if all
particles are also killed upon escaping $[v-KL, v+KL]$, then for all
sufficiently large $L$ (depending solely on $\phi$ and $\eps$), with
probability at least $(1-\eps) \theta$, the following hold:

\begin{longlist}
\item the erased $\mathit{KBRW}_1$ gets at least $f_0(L)$ descendants
in at most\break $c \log f_0(L)$ generations for some $c>0$,
\item $f_0(L)/K$ of them are in $\cN_v$.
\end{longlist}
\end{lemma}
\begin{pf}
Let $\tau$ be the first time that the killed branching random walk has
more than $2f_0(L)$ descendants. Let us show that the associated erased
branching random walk has at least $f_0(L)$ individuals at that point
with high probability. To see this, we first observe that by (i)
in\vadjust{\goodbreak}
Lemma \ref{Lkbrw1} the number of generations, $\tau$, is at most $c
\log f_0(L)$ for some $c>0$. Before time~$\tau$, for each new vertex
added to the branching random walk, the probability that it is a
self-intersection or hits an element of $F$ is no more than
$(2f_0(L)+|F|)/(2L-1)$. Thus the probability that a particular ray of
no more than $c \log f_0(L)$ generations contains a self-intersection
is, by Markov's inequality, at most
\[
\frac{c (\log f_0(L)) (2f_0(L)+|F|)}{2L-1} \to0
\]
as $L \to\infty$. Therefore, the number of vertices that are present
in the $\mathit{KBRW}_1$ but not in the erased $\mathit{KBRW}_1$ is, by Markov's
inequality again, at most $(1/2)\times (2f_0(L)) = f_0(L)$ with high
probability. We shall denote by $\mathit{EKBRW}_1$ the erased $\mathit{KBRW}_1$ which has
a set $F$ of forbidden locations.

By Lemma \ref{Lkbrw1}, we also know that $2f_0(L)/K$ individuals of
the $\mathit{KBRW}_1$ population are located in $\cN_v$. Since we have just
shown that the total number of individuals not in $\mathit{EKBRW}_1$ is
$o(f_0(L))$ with high probability, we deduce that at least $f_0(L)/K$
individuals of $\mathit{EKBRW}_1$ are located in $\cN_v$.
The proof of the lemma is complete.
\end{pf}

%s4.3 #&#
\subsection{Breadth-first search explorations}
The next three lemmas give us some information on the breadth-first
search exploration of a component $C_v$ of a given vertex $v \in V$ in
the random graph $G(t)$. For reasons that will soon become clear, we
wish to assume that by the point we start exploring the component
$C_v$, part of the graph has already been explored (a vertex has been
explored once all its neighbors have been observed). The part that has
already been explored (denoted~$F$) represents forbidden vertices, in
the sense that since we have already searched $F$, the breadth-first
search of $C_v$ can no longer connect to it.

We now specialize to the case where $\phi$ is the generating function
of a~$\operatorname{Poisson}(c)$ distribution, and in all that follows we let $\theta=
\theta(c)$ be the survival probability of a $\operatorname{Poisson}(c)$ Galton--Watson
tree. It turns out that we need to separately treat the case where $L$
is very close to $n$, and this will be done later in Lemma \ref{Llargepairs}.
%
%le9 #&#
\begin{lemma}
\label{Lsurvivalestin0} Fix $c>1$, $\eps>0$, $v\in V$ and fix $K =
K(\eps, c)$ as in Lemma \ref{Lsurvival-erased}.
We assume that a set $F$ containing at most $L/\omega(L)$ vertices
have already been discovered in $[v-KL,v+KL]$ (and $v$ is not one of them).
Then for all $n$ large enough (depending only on $\eps$ and $c$), if
$L< n / (2K)$, then with probability at least $\theta(1-\eps)$, a
search procedure of $C_v$ can uncover at least $f_0(L)/K$ vertices of
$C_v$ in $\cN_v$ without exploring more than $f_0(L)$ vertices in total
in $[v-KL,v+KL]$, and none outside.
\end{lemma}
\begin{pf}
Consider the breadth-first search exploration of $C_v$, with the
following modifications. We stop exploring the descendants of any
vertex outside of $[v-KL, v+KL]$. We\vadjust{\goodbreak} also completely stop the
exploration when more than $f_0(L)$ vertices have been discovered.
Also, we stop exploring the descendent of any vertex in $F$ and we fix
$d\ge1$ and truncate the offspring progeny at $d$, so that if an
individual has more than $d$ offspring, only the first $d$ encountered
are explored. This keeps the degree of any node in the genealogical
tree bounded by $d+1$. We choose $d=d_K$ as in Lemma \ref{Lkbrw1}. We
also stop exploring the descendants of an individual if the time
elapsed since the last time an ancestor of this individual visited $\cN
_v$ exceeds $\lambda=\lambda_{K,d}$, where $\lambda_{K,d}$ is as in
Lemma \ref{Lkbrw1}. We refer to this process as $\mathit{KBFS}_1$. More
formally, we use the following algorithm:\vspace*{8pt}

\textit{Step} 1. Set $\Omega_E=\varnothing$, $\Omega_A=\{v\}$. These
correspond to the explored and active vertices, respectively.

\textit{Step} 2. If $|\Omega_E|\ge f_0(L)$ we stop. Otherwise we
proceed to Step 3.

\textit{Step} 3. Set $\Omega_N=\varnothing$. For each $w\in\Omega_A$,
add its neighbors (excluding the parent of $w$) to $\Omega_N$ until $d$
have been added, or there are no more.

\textit{Step} 4. Add the vertices in $\Omega_A$ to $\Omega_E$.

\textit{Step} 5. Set $\Omega_A=\Omega_N\setminus\{\Omega_E\cup F\}$.
If $\Omega_A=\varnothing$, then we stop.

\textit{Step} 6. Remove from $\Omega_A$ all vertices outside of
$[v-KL,v+KL]$ and those that do not have an ancestor in $\cN_v$ fewer
than $\lambda$ generations away.

\textit{Step} 7. Go to Step 2.\vspace*{8pt}

This exploration can be exactly coupled with the $\mathit{EKBRW}_1$ considered
up to the first time $\tau$ that the total population size exceeds
$f_0(L)$, by taking in Lemma \ref{Lsurvival-erased} the set $F$ as it
is defined here. Lemma \ref{Lsurvivalestin0} thus follows directly
from Lemma \ref{Lsurvival-erased}.
\end{pf}

To establish the existence of a connection between two vertices $v$ and~$w$,
it will be useful to add another twist to the breadth-first search
exploration of~$C_v$ and~$C_w$, by \textit{reserving} some of the
vertices we discover along the way. That is, we decide not to reveal
their neighbors until a later stage, if necessary. This allows us to
keep a reserve of ``fresh'' vertices to explore at different locations
and that we know are already part of $C_v$ or $C_w$. To be more
precise, let $\eps>0$. Let $1<c'<c$ be such that $\theta(c') \ge
\theta
(c)(1-\eps/2)$. Let~$\nu$ be small enough that $c(1-\nu)>c'$. When
exploring $C_v$ through a method derived from breadth-first search, we
choose which vertices to reserve as follows: for each new vertex that
we explore,
if it has any offsprings, we choose one uniformly at random, and
reserve it with probability $\nu$ independently of anything else. (See
below for a rigorous formulation.) Note, in particular, that the set of
vertices that get reserved is dominated by a Poisson thinning of the
original exploration procedure, with thinning probability $\nu$.
Let $K=K(\eps/2,c')$ be as in Lemma \ref{Lsurvival-erased}.
Note that with this choice of $\nu$ and $K$, the survival probability~$\theta'$ of $\mathit{EKBRW}_1$ is at least
%
%e13 #&#
\begin{equation}\label{reserve}
\theta'(c) \ge\theta(c')(1-\eps/2) \ge\theta(c)(1-\eps)
\end{equation}
for all $L$ sufficiently large (depending solely on $c$ and\vadjust{\goodbreak}
$\eps$).

Thus, starting from a vertex $v$, a branching random walk killed when
escaping $[v-KL, v+KL]$ with this reservation procedure survives
forever with probability at least $\theta(c)(1-\eps)$. From this we
deduce without too much trouble the following result.
%
%le10 #&#
\begin{lemma}\label{reserving}
Fix $c>1$, $\eps>0$, $v\in V$. Let $\nu= \nu(\eps)$ as above, and
assume that a set $F$ containing no more than $L/ \omega(L)$ vertices
have been discovered.
Then for all sufficiently large $n$ (depending solely on $\eps$ and
$c$), if $L<n/(2K)$, the following hold with probability at least
$(1-2\eps)\theta$:
\begin{longlist}
\item A search procedure can uncover at least $f_1(L)=
f_0(L)/K$ vertices in~$\cN_v$ without uncovering more than $2f_0(L)$
vertices in total in $[v-KL,\allowbreak v+KL]$, and none outside.

\item At least $\delta f_0(L)$ vertices are reserved in $\cN
_v$, for $\delta= (1-e^{-c})\nu(\eps)/(4K)$.
\end{longlist}
\end{lemma}
\begin{pf}
We apply the above reservation method to $\mathit{KBFS}_1$ (see the proof of
Lemma \ref{Lsurvivalestin0}). Formally, we introduce a set $\Omega_R$
of reserved vertices (initially $\Omega_R=\varnothing$). We use the
same algorithm as for the modified breadth-first search but now Step 7
becomes:\vspace*{8pt}

\textit{Step} 7$'$. Partition $\Omega_A$ into classes of vertices with
the same parent in the exploration. Choose uniformly from each class a
representative and with probability $\nu$ this representative is added
to $\Omega_R$ and removed from $\Omega_A$. Go to Step 2.\vspace*{8pt}

We call this new search procedure $\mathit{KBFS}_2$. Let $\tau$ be the time we
have discovered $f_0(L)$ nonreserved vertices. At this time the total
number of explored vertices is less than $2f_0(L)$ and thus, similar to
the proof of Lemma~\ref{Lsurvivalestin0}, we can couple the
exploration with an erased $\mathit{KBRW}_1$ where the offspring distribution
has a slightly modified offspring distribution (a randomly chosen
offspring is removed with probability $\nu$). We call this an erased $\mathit{KBRW}_2$.
Reasoning as in Lemma \ref{Lsurvivalestin0}, and using (\ref
{reserve}), we see that (i) holds with probability at least $\theta
(1-\eps)$, provided that $n$ is large enough and $L< n/(2K)$.
For each new vertex exposed by $\mathit{KBFS}_2$ in $\cN_v$, it has a reserved
offspring in~$\cN_v$ with probability at least $(1-e^{-c})\nu/2$, as
if $u \in\cN_v$ and $X$ are uniformly distributed on $\{-L, \ldots,
-1, 1, \ldots, L\}$, then $u+X \in\cN_v$ with probability at least
$1/2$. Thus (ii) follows from (i) and from Chebyshev's inequality.
\end{pf}

With this lemma we are now able to show that a vertex $v$ connects to
$v\pm KL$ with probability essentially $\theta$, and that many vertices
in the same component may be found without revealing too much inside
$[v-KL,\break v+KL]$.
%
%le11 #&#
\begin{lemma}\label{Lsurvivalestout}
Fix $c>1$, $\eps>0$, and let $K$ be as in Lemma
\ref{Lsurvivalestin0}. Let $0<\zeta<\zeta'<1/2$ and let $v \in V$. Assume
that a set $F$ of no more than $L/\omega(L)$ vertices have already
been\vadjust{\goodbreak}
explored in $[v-KL, v+KL]$ and $v$ is not one of them. Let~$\cB_{K,v}$
denote the event that $v$ is connected to at least $L^{\zeta}$
unexplored vertices in the range $[v+KL, v+(K+1)L]$ which may be
discovered by searching no more than $L^{\zeta'}$ vertices. Then for
all sufficiently large $n$ (depending solely on $\eps$, $c$, $\zeta$
and $\zeta'$), if $L<n/(2(K_0(\eps)+1))$, then
\[
\P(\cB_{K,v}) \ge\theta(c)(1-\eps).
\]
\end{lemma}
\begin{pf}
Consider the $\mathit{KBFS}_2$ exploration of $C_v$, stopped when a total of
$f_0(L)=L^{\zeta'}/2$ vertices of $C_v$ have been exposed (additional
to those exposed initially).
%We will not see any self-intersections or connections with the set of
%already exposed vertices with probability at least $1-L^{-1/4}/4$, and
%thus with probability tending to 1 we can couple this exploration
%exactly with the reserved branching random walk, considered until the
%time that the total progeny exceeds $f_0(L)$.
By Lemma \ref{reserving}, with probability at least $\theta(c)(1-\eps)$
if $n$ is large enough and $L< n/ (2K)$, this search reveals at least
$k=\delta f_0(L)$ reserved vertices within~$\cN_v$, and no more than
$L^{\zeta'}$ vertices in the range $[v-KL,\break v+KL]$ have been explored
(let $\cA_{v}$ denote this event).
On $\cA_{v}$, label $v_1, \ldots, v_k$ the first $k$ such vertices to
have been discovered in $\cN_v$. After this stage, we then continue the
$\mathit{KBFS}_2$ exploration started from $\{v_1, \ldots, v_k\}$ only, until a
total of $L^{\zeta'}/2$ further vertices are exposed. Note that the
exploration can be coupled with a system of $k$ erased $\mathit{KBRW}_2$,
started from $v_1, \ldots, v_k$. The total number of vertices searched
by the end of the second stage will be no more than $f_0(L) + L^{\zeta
'}/2 \le L^{1/2}$. Thus, as in Lemma~\ref{Lsurvival-erased}, the
probability each\vspace*{1pt} particular vertex gives rise to a self-intersection is
no more than $(|F|+\sqrt L)/(2L-1)\to0$ as $n\to\infty$.

Moreover, using domination by a branching process (Lemma \ref
{Lbreadth}), it is easy to see that the number of generations for the
$L^{\zeta'}/2$ vertices to be discovered by the branching random walk
is at least $b\log L$ for some $b>0$, with probability tending to 1 as
$n \to\infty$. Now, for every $1\le i \le k$, the probability that the
erased branching random walk (erased $\mathit{KBRW}_2$) started from $v_i$ has
a~descendant that hits $[v+KL, v+(K+1)L]$ in fewer than $b \log L$ steps
is at least $\theta(c)(1-o(1)) \ge\theta(c)/2$ for $n$ large enough
(depending solely on $\eps$ and~$c$). Thus we deduce that, on the event
$\cA_{v}$, the number of particles that hit $[v+KL,v+(K+1)L]$
stochastically dominates a binomial random variable
$\operatorname{Bin}(\delta f_0(L),
\theta(c)/2)$. By applying Chebyshev's inequality, this is with high
probability greater than $\delta f_0(L) \theta(c)/4 \ge L^{\zeta}$ for
all $n$ sufficiently large (depending on $c, \eps, \zeta, \zeta'$).
When this occurs, $\cB_{K,v}$ holds, so the result follows.
\end{pf}

Let $c>1$, $\eps>0$ and fix $K=K_0(c,\eps)$ as in Lemma \ref
{Lsurvivalestout}. We now prove that if the connected component of a
given vertex is not finite then it must spread more or less uniformly
over $V$. As desired, this is achieved by keeping a density of explored
sites small, lower than $1/\omega(L)$. Let $v \in V$ and split the
vertex set $V$ into $r+1=\lceil n/(KL)\rceil$ disjoint strips $(I_0,
\ldots, I_r)$ of size~$KL$, except for the last one which may be of
size smaller than $KL$. Let $J_i$ denote the initial segment of $I_i$
of length $L$. Since $L(n) \ge(\log n)^{2+\xi}$ for some positive~$\xi
$ by assumption, we may find $\zeta<1/2$ such that $L(n) >(\frac
{4}{\theta}\log n)^{1/\zeta}$. Let $\zeta' = (\zeta+ 1/2)/2$, hence
$\zeta<\zeta' < 1/2$.
%
%le12 #&#
\begin{lemma}\label{Lsurvivalcrossing} With the above notation, assume
that no more than $L^{2\zeta'}$ vertices have already been exposed in
each strip $I_0,\ldots,I_r$. Let $\cC_{K,v}$ denote the event that $v$
is connected to at least $k= L^{\zeta}$ vertices in strips $I_3,
\ldots
, I_{r-1}$, which may be discovered without exposing more than an
additional $L^{2\zeta'}$ vertices in each strip, and that in each $J_i,
3\le i \le r-1$, at least $k/2$ vertices connected to $v$ are
unexplored at the end of the search procedure. Then
\[
\P(\cC_{K,v} ) \ge\theta(c)(1-\eps)
\]
for all $n$ large enough (depending solely on $\eps$ and $c$), and
provided $L<n/(2K)$.
%Moreover,
%the same result holds with reserved breadth-first search with
%reservation probability $\nu(c,\eps)$.
\end{lemma}
\begin{pf}
This is basically proved by iterating Lemma \ref{Lsurvivalestout}. We
can assume that $v$ is in strip $I_1$. In the first step we explore
$C_v$ using $\mathit{KBFS}_2$ killing at the boundary of $I_0\cup
I_1\cup I_2$. The arguments of Lemma \ref{Lsurvivalestout} still carry
through to obtain that $\cB^0_{K,v}$ holds with probability at least
$\theta (c)(1-\eps)$ where $\cB^0_{K,v}$ is the event that $v$ is
connected to at least $k$ unexplored vertices in $J_3$, and fewer than
$L^{\zeta'}$ vertices are explored in finding them.

Then define inductively for $1\le i \le r-4$, $\cB^i_{K,v} = \cB
^{i-1}_{K,v} \cap\cC^i_{K,v}$, where $\cC^i_{K,v}$ is defined as
follows. On $\cB^{i-1}_{K,v}$, let $v_1, \ldots, v_k$ be a list of $k$
vertices in the range $J_{i+2}$ that are the first to be discovered in
this search procedure in this range. Then $\cC^i_{K,v}$ is the event
that we can find at least $k$ connections between $v_1, \ldots, v_k$
and $J_{i+3}$ without exploring more than an additional $L^{2\zeta'}$
new vertices in $I_{i+2}\cup I_{i+3}$.

This is where\vspace*{2pt} it starts to pay off to allow for the exploration to
unfold in a~partially revealed environment in Lemma
\ref{Lsurvivalestout}. Indeed, let $i\ge1$ and condition on
$\cB^{i-1}_{K,v}$. We reserve (i.e., do not explore further)
$v_{{k/2+1}}, \ldots, v_k$. We explore successively the components
$C_{v_1}, \ldots, C_{v_{k/2}}$, each time performing $\mathit{KBFS}_2$
of Lemma~\ref{Lsurvivalestout}. Since we never reveal more than
$L^{\zeta'}$ vertices at each of those $k/2$ steps, and since we did
not reveal more than $(k/2)L^{\zeta'} =
L^{\zeta+\zeta'}/2<L^{2\zeta'}/2$ other vertices in $I_{i+2}$
previously (since $\cC^{i-1}_v$ holds), we see that the search may be
coupled with high probability (depending solely on $c$ and $\eps$) to
$(k/2)$ erased $\mathit{KBRW}_2$ started at $v_1, \ldots, v_{k/2}$.
Thus the total number of connections between $v_1, \ldots, v_{k/2}$, to
$J_{i+3}$ is dominated from below by $k \operatorname{Bin}(k/2,
\theta(c)/2)$. Indeed, for each of $(k/2)$ trials there is a
probability $\theta(c)(1-\eps) \ge\theta(c)/2$ of success (by Lemma
\ref{Lsurvivalestout}), in which case $k$ connections are added. Thus,
using standard Chernoff bounds on binomial random variables,
\begin{eqnarray*}
\P((\cB^i_{K,v})^\complement| \cB^{i-1}_{K,v}) &=& \P((\cC
^i_{K,v})^\complement| \cB^{i-1}_{K,v}) \\
& \le &\P\bigl( k \cdot\operatorname{Binomial}\bigl(k/2, \theta(c)/2\bigr) < k\bigr) \\
&\le&\exp\bigl( - \theta(c) k/4\bigr).
\end{eqnarray*}
It follows by easy induction that for all $n$ large enough, letting
$\cC
'_{K,v} = \bigcap_{i=0}^r \cB^i_{K,v}$,
\[
\P(\cC'_{K,v}) \ge(1- e^{-\theta k/4})^r\P(\cB^0_{K,v}).
\]
Since $L> ( (4/\theta)\log n)^{1/\zeta}$, $r = \lceil n/(2LK)\rceil$
and $k=L^\zeta$, it follows that $(1- e^{-\theta k/4})^r \sim\exp(-r
e^{-\theta k/4}) \ge\exp(-r/n) \to1$. Vertices can only be discovered
during two consecutive steps of the proof, and hence, the total number
of vertices discovered in each strip is no more than $L^{2\zeta'}$.
Thus $\cC_{K,v} \supset\cC_{K,v}'$. The proof of the lemma is complete.
\end{pf}
%
%le13 #&#
\begin{lemma}
Let $\cD_v =\{ |C_v|> \log L\}$. Then for any fixed $v \in V$,
\[
\P(\cD_v) \to\theta(c)
\]
and for $v,w$ fixed in $V$,
\[
\P(\cD_v \cap\cD_w) \to\theta(c)^2.
\]
\end{lemma}
\begin{pf}
This is a direct consequence of Lemma \ref{Lbrwcoupl} and the remark
following it.
\end{pf}
%
%le14 #&#
\begin{lemma}\label{Lpairs}
%Fix $0<\zeta<\zeta'<1/2$ and assume $L > ((4/\theta)\log n)^{1/\zeta}$.
Fix $c>1$, $\eps>0$. Let $v,w$ be chosen uniformly at random in $V$,
and let $\cE= \{C_v = C_w\}$ be the event that they are connected.
If $L< n/(5K)$, and $n$ is large enough (depending solely on $\eps$
and $c$) then
%
%e14 #&#
\begin{equation}\label{LBconnec3}
\P( \cE^\complement| \cD_v \cap\cD_w) \le\eps.
\end{equation}
\end{lemma}
\begin{pf}
We fix $K(c,\eps)$ as in Lemma \ref{Lsurvivalcrossing}. We apply this
lemma a first time by exploring $C_v$ as specified in this lemma with a
set of forbidden vertices (vertices previously explored) being empty.
We then let $F$ be the set of all vertices explored during that procedure.

We apply one more time Lemma \ref{Lsurvivalcrossing} by exploring
$C_w$ using a set of forbidden vertices given by $F$ (which must
necessarily satisfy the assumptions of Lemma~\ref{Lsurvivalcrossing},
since the search of $C_v$ did not reveal more than $L^{2\zeta'}$
vertices in each strip). Note that conditionally given $\cD_v \cap\cD
_w$, both $\cC_{K,v}$ and $\cC_{K,w}$ must hold with high probability
(depending solely on $c$).
Let us show that~$\cE$ must then hold with high probability.

Since $\cC_{K,v}$ and $\cC_{K,w}$ hold, we know that each interval
\mbox{$J_i, 3\le i \le r-1$}, contains at least $L^\zeta/2$ unexplored
vertices from both $C_v$ and $C_w$. We now apply Lemma \ref{reserving}
repeatedly, starting from each of these unexplored vertices. Since
$L<n/5K$, we have that $r \ge4$. While fewer than $\delta f_0(L)$
vertices have been reserved in $J_i$, we know that fewer than $f_0(L)$
vertices have in total been explored and thus Lem\-ma~\ref{reserving}
can still be applied. We deduce that (conditionally given
$\cD_v\cap\cD_w$) in each $J_i, 3\le i \le r-1$, with probability
greater than $1-o(1)$ depending solely on $\eps$ and~$c$, there are
$\delta L/\omega(L)$ reserved\vadjust{\goodbreak} vertices from $C_v$ and $\delta
L/\omega(L)$ reserved vertices from $C_w$, with $\delta$ as in
Lemma~\ref{reserving}. Thus at least one $J_i$ (say $J_1$) contains $\delta
L/\omega(L)$ unexplored vertices from both $C_v$ and $C_w$. The
probability to not observe a connection between these $(\delta L
/\omega (L))^2$ pairs of vertices inside $J_1$ is at most (by revealing
only the status of the edges connecting each such pair)
\begin{eqnarray*}
\biggl(1- \frac{c}{2L}\biggr)^{(\delta L /\omega(L))^2} &\le&\exp
\biggl(-
\frac{c}{2L} \frac{\delta L^2}{\omega(L)^2}\biggr)\\
& = &\exp\biggl( - \frac{c \delta L}{2 \omega(L)^2} \biggr) \to0
\end{eqnarray*}
for all $n$ sufficiently large. Thus $C_v= C_w$ with high probability
(depending solely on $\eps$ and $c$) given $\cD_v\cap\cD_w$, and hence,
$\cE$ holds with high probability (depending solely on $\eps$ and $c$)
given $\cD_v\cap\cD_w$.
\end{pf}

We deal with the case $L \ge n/(5K)$ separately.
%
%le15 #&#
\begin{lemma}\label{Llargepairs}
Fix $c>1$, $\eps>0$. Let $v,w$ be chosen uniformly at random in $V$,
and let $\cE= \{ C_{v} = C_{w}\}$ be the event that they are
connected. If $L\ge n/(5K)$, and $n$ is large enough (depending
solely on $\eps$ and $c$) then
%
%e15 #&#
\begin{equation}\label{LBconnec3large}
\P( \cE^\complement; \cD_v \cap\cD_w) \le\eps.
\end{equation}
\end{lemma}
\begin{pf}
We let $\bar\cE$ denote the event that we can explore $C_v$ until at
most~$L^{0.7}$ vertices have been exposed finding at least $L^{0.6}$
reserved vertices and also can explore $C_w$ until at most $L^{0.45}$
vertices have been exposed finding at least $k=L^{0.44}$ reserved
vertices. By a simple modification of Lemma~\ref{reserving}, the
probability of this event given $\cD_v \cap\cD_w$ is at least $1-o(1)$.

Let us show that, given the above event $\bar\cE$, $C_{v}$ and $C_{w}$
intersect with high probability. We partition $\{1,\ldots,n\}$ into
$s=5 K +1 $ disjoint intervals of size less than or equal to $L$. On
the above event, by the pigeonhole principle there must be at least one
region of size at most $L$, denoted $I$, with more than~$L^{0.6}/s$
reserved vertices from $C_v$.
We denote by $w_1, \ldots, w_k$ the $k$ reserved vertices from $C_w$.
For each $1\le i\le k$, we continue to explore $C_{w_i}$ by
breadth-first search for $6(s+1)^2$ generations, or until a descendent
is observed in interval $I$. Since $s$ depends only on $\eps$ and $c$,
with probability at least $1-o(1)$ depending on $\eps$ and $c$, no
self-intersections occur throughout this evolution.
We claim that the probability the evolution of $C_{w_i}$ results in us
finding a~descendent in $I$ is at least $\theta(c)/3$ for $n$ large
enough depending solely on $\eps$ and $c$. Indeed this occurs if we can
find a ray emanating from $w_i$ where the corresponding random walk
goes around the circle in less than $6(s+1)^2$ levels. We let
$(X_j)_{j\ge1}$ denote the location of a random walk on $\mathbb{Z}$
which starts at~0 and where the jump distribution is uniform on $\{
-L,\ldots,-1,1,\ldots,-L\}$. It is clear that $(X_j)_{j\ge0}$ is a
martingale, as is
\[
\bigl(X_j^2-jL(L+1)(2L+1)/(6L)\bigr)_{j\ge0}.\vadjust{\goodbreak}
\]
Letting $T$ denote the time the walk goes above $sL$, or below $-sL$,
we see that by optional stopping $\bbe(T)<3(s+1)^2$. Thus, by Markov's
inequality, $\bbp(T>6(s+1)^2)<1/2$. Hence, a random walk on $V$ with
jumps uniformly distributed on $\{-L,\ldots,-1,1,\ldots,-L\}$ goes
around the circle in less than $6(s+1)^2$ steps with probability at
least $1/2$.
It follows that the desired ray exists with probability at least
$(1-o(1))\theta(c)/2$, depending solely on $\eps$ and $c$. Thus, given
$\bar\cE$, we can find $L^{0.43}$ reserved vertices from $C_w$ and
$L^{0.6}$ reserved vertices from $C_v$ in the interval~$I$, with
probability $1-o(1)$ depending solely on $\eps$ and $c$. Looking one
level further, the number of connections between~$C_v$ and $C_w$ in
this region is $\operatorname{Bin}(L^{1.03}$, $c/(2L))$ which is
larger than 1
with probability $1-o(1)$ depending solely on $\eps$ and~$c$. Equation
(\ref{LBconnec3large}) now follows.
\end{pf}

We are now ready to finish the proof of (ii) in Theorem \ref{TuniqGC}.
\begin{pf*}{Proof of Theorem \ref{TuniqGC}\textup{(ii)}}
Let $v, w$ be chosen uniformly on $V$. Let $\cE= \{ C_{v} = C_{w}\}$
be the event that they are connected.

Consider $W = \{v \in V\dvtx \cD_{v} \mbox{ holds}\}$. We already know that
$\E(|W|/n) \to\theta$ and $\E(|W|^2/n^2) \to\theta^2$, so that
\[
\frac{|W|}n \to\theta
\]
in probability. Furthermore, observe that if $v,w$ are uniformly chosen
in $W$, then $C_v = C_w$ with high probability depending solely on $c$
by Lemmas \ref{Lpairs} and \ref{Llargepairs}. Also, if $v \in W$,
then clearly $C_v \subset W$. Hence, $W$ consists of a union of
clusters. Let $X_n$ denote the size of a cluster from $W$ chosen
according to size-biased picking, that is, $X_n =_d |C_v|/ |W|$, where
$v$ is chosen uniformly at random in $W$.
It is a well-known consequence of exchangeability (and easy to see) that
\[
\P(C_v = C_w |v,w \in W) = \E(X_n).
\]
By (\ref{LBconnec3}), if $L< n/ (4LK)$ [resp., by (\ref
{LBconnec3large}) if $L \ge n/ (4LK)$], we have that
\[
\P(C_v = C_w | v,w \in W) \to1,
\]
hence, $\E(X_n) \to1$. Since $X_n \le1$, it follows that $X_n \to1$
in probability. This implies that, for all $\eps>0$, with high
probability depending solely on $\eps$ and~$c$, $W$~contains a
component of size at least $\frac{|W|}{n}(1-\eps) \ge\theta(c)(1-2
\eps
)$. This proves the existence of a giant component of mass relative to
$V$ equal to $\theta$ in the limit $n \to\infty$.

Let us show that all other components are small. Note that by the
above, we already know that the second largest component size, $L^2_W$,
is such that $L^2_W / |W| \to0$ in probability. Hence,\vspace*{-1pt}
$L^2_W / n \to 0$ as well. Let $L^1_{W^\complement}$ be the largest
component size in $W^\complement$. By definition,\vadjust{\goodbreak} $L^1_{W^\complement}$
is smaller in size than $\log L$. Since $L^2_n
\le\max(L^1_{W^\complement}, L^2_W)$, we conclude that
\[
\frac{L_n^2}n \to0
\]
in probability, as desired. The proof of (ii) in Theorem \ref{TuniqGC}
is complete.\vspace*{-2pt}
\end{pf*}

We now conclude with the proof of (iii) in Theorem \ref{TuniqGC}.\vspace*{-2pt}
\begin{pf*}{Proof of Theorem \ref{TuniqGC}\textup{(iii)}}
Since $L=o(\log n)$, we have $L<\frac{4a}{6c}\log n$, with $a>0$ as in
the statement of (iii). We begin by dividing $1,\ldots,n$ into
$n^{1-a}\log n$ disjoint intervals of size $n^a/\log n$, labeled
$A_1,\ldots, A_{n^{1-a}\log n}$. In each interval we show that we can
find an interval of size $L$, none of whose vertices have been involved
in a transposition by time $t$ with high probability. We show in fact,
that all the $n^{1-a}\log n$ intervals contain such a sub-interval with
high probability. Thus the largest component must be of size smaller
than $2n^a/\log n$, and hence, in particular, there will be no giant components.

For a given interval of size $L$, the number of potential edges
connected to vertices in this interval is $2L^2-
{L\choose2}=(3L^2+L)/2$. Each\vspace*{1pt} of these edges is present with probability
$c/(2L)$. We call the interval empty if none of the edges are present.
The probability a given interval of size $L$ is empty is
\[
\bigl(1-c/(2L)\bigr)^{(3L^2+L)/2}\sim\exp\biggl(-\frac{c}{4}(3L+1)\biggr).
\]
We divide each $A_i$ into $\lfloor n^a/(L\log n)\rfloor$ intervals of
size $L$, denoted
\[
A_i^1,  A_i^2,  \ldots,  A_i^{\lfloor n^a/(L\log n)\rfloor}.
\]

Let $S_i^{2k}=\{A_i^{2k}\mbox{ is empty}\}$. We consider the set of
events
\[
\bigl\{S_i^{2k},   1\le k\le\tfrac{1}{2}\lfloor n^a/(L\log n)\rfloor\bigr\},
\]
which are independent since each interval is at distance at least $L$
from any other.
For each $i$, we let
\[
B_i=\sum_{k=1}^{\lfloor n^a/(L\log n)\rfloor/2}\mathbf
{1}_{\{S_i^{2k}\}}
.
\]
We have
\[
\P(B_i>0)\sim1-\exp\biggl[-\frac{1}{2}\biggl\lfloor\frac
{n^a}{L\log
n}\biggr\rfloor e^{-c(3L+1)/4}\biggr]
\]
and so
\begin{eqnarray*}
&&\P(B_i>0\mbox{ for all }1\le i\le n^{1-a}\log n)\\[-2pt]
&&\qquad \sim\exp\biggl[-n^{1-a}\log n\exp\biggl(-\frac{1}{2}
\biggl\lfloor
\frac{n^a}{L\log n}\biggr\rfloor e^{-c(3L+1)/4}\biggr)\biggr]\\[-2pt]
&&\qquad \to1\qquad\mbox{as } n\to\infty
\end{eqnarray*}
since $L<\frac{4a}{6c}\log n$. The proof of Theorem \ref{TuniqGC} is
complete.\vadjust{\goodbreak}
\end{pf*}

%s5 #&#
\section{\texorpdfstring{Proof of Theorem \protect\ref{Tdist}}{Proof of Theorem 2}}
\label{Sproofofbio}

We will prove a stronger result than Theorem \ref{Tdist} by allowing
the distribution of edge-lengths to be more general. Recall the
definitions of $(p_\ell)$ and $\eps_n$ given at the beginning of
Section \ref{Scomponents}. Let $(\tau_i)_{i \ge1}$ be a sequence of
i.i.d. transpositions with $\tau_1 = (i \ j)$ where $i,j$ are chosen
uniformly from $\{u,v \in V\dvtx \|u-v\|=L\}$. Then we construct the
permutation $\sigma_t = \tau_1 \circ\cdots\circ\tau_{N_t}$, where
$(N_t, t \ge0)$ is an independent Poisson process. In words, at rate
1, we transpose two markers at random with distance $D$, where $D$ is
chosen according to the distribution $(p_\ell)$. We recover the process
$(\sigma_t\ge0)$ when $(p_\ell)$ is the uniform distribution on $\{1,
\ldots, L(n)\}$.\vspace*{-2pt}
%
%th6 #&#
\begin{theorem}\label{Tdistgeneral}
Assume $\eps_n\to0$ as $n\to\infty$. Then we have the following
convergence in probability as $n\to\infty$: for all $c>0$,
\[
\frac1{n}\delta(\sigma_{cn/2})\to u(c).\vspace*{-2pt}
\]
\end{theorem}

There is a natural coupling between the process $(\sigma_t,t\ge0)$ and
the random graph $(G(t),t\ge0)$ defined in Section \ref{Scomponents}.
The coupling is an adaptation of the coupling with the Erd\H{o}s--Renyi
random graph in \citet{beres}.
%We also construct a multigraph denoted by $\bar G(t)$.
Consider the following procedure. Initially, $G_0$ consists of isolated
vertices.
Suppose that at time $t$, a transposition $\tau=(i,j)$ is performed.
If~$G(t)$ already contains the edge $(i,j)$ we do nothing, else we add it
to the graph.

The relationship between $\sigma_t $ and $G(t)$ is not one-to-one;
however, the following deterministic observation holds as in
\citet{beres}. For every \mbox{$t\ge0$}, every cycle of $\sigma_t$ is a subset of a
certain connected component of $G(t)$. That is, the partition of $V$
obtained from considering the cycle decomposition of $\sigma_t$ is a
refinement of the partition obtained from considering the connected
components of~$G(t)$. This is easily proved by induction on the number
$N_t$ of transpositions up to time $t$, after observing that the cycle
decomposition of $\sigma_t$ undergoes a coagulation-fragmentation
process. Indeed, every transposition $(i,j)$ that involves two
particles from the same cycle yields a fragmentation of that cycle,
while if the two particles are in distinct cycles they merge.

This coupling is the basis of our proof. Armed with Lemmas \ref
{LNclusterRG} and \ref{LconcX}, in order to prove Theorem \ref
{Tdistgeneral} we need to show that $K_t$ and $|\sigma_t|$ differ
by~$o(n)$ (Lemma \ref{Lcycle-rg}), where we recall that $|\sigma|$ is the
number of cycles of the permutation $\sigma$.\vspace*{-2pt}
%
%le16 #&#
\begin{lemma}\label{Lcycle-rg}
%Let $K_t$ and $|\sigma(t)|$ be the number of components and the number
%of cycles at time $t$, respectively.
Assume $\eps_n\to0$. Let $t=cn/2$, where $c>0$. As $n\to\infty$,
\[
\frac{|\sigma_t|-K_{t}}{n}\to0
\]
in probability.\vspace*{-2pt}
\end{lemma}
\begin{pf}
This argument is somewhat analogous to the proof of Lemma~6 in
\citet{nb}. First we note that by the properties\vadjust{\goodbreak} of the coupling between
$\sigma_t$ and $G(t)$, it is with probability 1 the case that $K_t \le
|\sigma_t|$.
To prove a bound in the converse direction, we need to distinguish
between small and large cycles or components. We say that a cycle of
$\sigma_t$ or a component of~$G(t)$ is \textit{small} if it has a size
less than $1/\sqrt\eps_n$ and \textit{large} if it has size at least
$1/\sqrt\eps_n$.

Note that the number of large cycles and the number of large components
is at most $n\sqrt\eps_n = o(n)$.
It thus suffices to control the difference between the number of small
cycles and the number of small components. However, note that at any
time, the probability of generating a small cycle by fragmentation is
at most $4 (1/\sqrt\eps_n)\eps_n$. To see where this comes from,
suppose the current permutation is $\sigma_t = \sigma$, and the first
position for the transposition $(i,j)$ to be performed has been chosen.
Thus $j$ will be one of the $n-1$ other vertices chosen according to
the distribution $(p_{\ell})$. Then to produce a cycle of size exactly
$k$, $j$ must be equal to $\sigma^{k-1}(i)$ or $\sigma^{-k+1}(i)$.
(Depending on the exact size of the cycle containing $i$, there may be
two other points allowed.) Thus, conditioning on the point $i$, the
probability of creating a fragment of size smaller than $1/\sqrt{\eps_
n}$ is at most $4 (1/\sqrt\eps_n)\eps_n$, as claimed. It follows that
since each excess small cycle must have been generated by such a
fragmentation at some time $s\le t$, and since transpositions occur at
rate 1,
\[
\bbe[|\sigma_t|-K_t]\le n\sqrt\eps_n + 4t\sqrt\eps_n.
\]
Thus by Markov's inequality, taking $t = cn/2$, for all $\delta>0$,
\begin{eqnarray*}
\bbp\biggl(\frac{|\sigma_t|-K_{t}}{n}>\delta
\biggr)&\leq&\frac{\bbe(|\sigma_t|-K_{t})}{\delta n}\\
&\leq&\sqrt\eps_n\frac{1+2c}{\delta} \to0
\end{eqnarray*}
as $n \to\infty$. The proof is complete.
\end{pf}
\begin{pf*}{Proof of Theorem \ref{Tdistgeneral}}
$\!\!\!$The proof of Theorem \ref{Tdistgeneral} now follows directly~from
Lemmas \ref{LNclusterRG}, \ref{LconcX} and \ref{Lcycle-rg}. Indeed,
$\delta(\sigma_t) = n- |\sigma_t|$. By Lemma \ref{Lcycle-rg},
\mbox{$n^{-1}(|\sigma_t| - K_t)$} tends to $0$ in probability. We have
concentration of $K_t$ around its mean by Lem\-ma~\ref{LconcX}, and the
mean is obtained in Lemma \ref{LNclusterRG}. Putting these pieces
together we obtain Theorem~\ref{Tdistgeneral}.\vspace*{-2pt}
\end{pf*}

%s6 #&#
\section{\texorpdfstring{Proof of Theorem \protect\ref{Tmicro}}{Proof of Theorem 3}}
\label{Sproofofmicro}

We consider the case where $L$ is bounded (say by some constant $C$)
and show that if $t=cn/2$ with $c>0$, then
$\delta(t)/n$ is bounded away from $c/2$, where we write $\delta(t) =
\delta(\sigma_t)$.
%
%le17 #&#
\begin{lemma}\label{Lboundedfrags}
Assume $L$ is bounded. Fix $c>0$ and let $t=cn/2$. Then there exists
$\eta= \eta_c>0$ such that $\delta(t)\le(1-\eta)cn/2$ with high
probability.
\end{lemma}

In the statement above and in what follows, the expression \textit{with
high probability} means with probability tending to 1 as
$n\to\infty$.
\begin{pf*}{Proof of Lemma \ref{Lboundedfrags}}
Since each transposition decreases the number of cycles by 1 if there
is a coagulation and increases it by 1 if there is a~fragmentation, we have
%
%e16 #&#
\begin{equation}\label{deltafrag}
\delta(t) = N_t - 2F_t,
\end{equation}
where $F_t$ is the total number of fragmentations by time $t$. It
suffices to show that $F_t \ge\eta n$ when $t=cn/2$, for some $\eta>0$.
Let $(i,j) \in\cR_L$.
Consider the event $\cA_{ij}$ that the transposition $(i,j)$ occurred
twice by time $t$, and that no other transposition involved either $i$
or $j$ by time $t$. There are $4L-2$ possible transpositions involving
$i$ or $j$ but not both, with each occurring at rate $1/(nL)$.
Thus the number of such transpositions that occur by this time is
Poi($t(4L-2)/nL$) which has a positive probability, $q_c$, of being 0.
Further, the number of times transposition $(i,j)$ occurs by time $t$
is Poi($t/nL$) and thus we have a positive probability, $p_c$, of it
occurring exactly twice. Thus $\P(\cA_{i,j}) =q_c p_c >0$ for each
$(i,j) \in\cR_L$.

Moreover, the events $(\cA_{2iL, 2iL+1})_{0\le i \le\lfloor
n/2L\rfloor
-1}$ are independent and each occurs with probability $q_cp_c$. Note
that the number $F_t$ of fragmentations satisfies
\[
F_t \ge\sum_{i=0}^{\lfloor n/2L\rfloor-1} \mathbf{1}_{\cA_{2iL, 2iL+1}}.
\]
It thus follows from Chebyshev's inequality that $\P
(F_t>nq_cp_c/(4C))\to1$, where $C$ is an upper-bound on $L$. Hence,
$F_t \ge\eta_c n$ with $\eta_c = q_c p_c/(4C)$. Plugging back in
(\ref{deltafrag}) completes the proof.
\end{pf*}

We now turn toward the proof of Theorem \ref{Tmicro}. Assume without
loss of generality that $L(n)=L$ is constant.
As $N_t=_d \operatorname{Poisson}(t)$, we obtain directly from
Chebyshev's inequality that $\frac{1}{n}N_{cn/2} \to c/2$ in probability.

It thus suffices to show that $\frac{1}{n}F_{cn/2}$ also has a limit as
$n \to\infty$. Let $(\cF_t)_{t\ge0}$ be the filtration associated
with the entire history of the process; that is, $\cF_t = \sigma(\tau
_i, i \le N_t)$. For $s\ge0$, let $g_n(s)$ denote the $\cF
_s$-measurable random variable giving the instantaneous rate of
fragmentation given $\sigma_s$. Let\looseness=1
\[
A_t = \int_0^t g_n(s)\,ds
\]\looseness=0
and observe that if $M_t = F_t - A_t$, then $(M_t, t\ge0)$ is a
martingale with respect to the filtration $(\cF_t)_{t\ge0}$, for each
$n \ge1$.

We prove convergence of $n^{-1}F_{t}$ (with $t=cn/2$) in two steps:
\begin{longlist}
\item $n^{-1} A_{t}$ converges,
\item $n^{-1} M_t \to0$ in probability, which will follow from
Doob's inequality.
\end{longlist}

Note first that by a change of variable,
%
%e17 #&#
\begin{equation}\label{Echgvar}
\frac1n A_{cn/2} = \frac12 \int_0^c g_n(sn/2) \,ds.
\end{equation}

%le18 #&#
\begin{lemma}
\label{LEA} There exists a nonrandom function $g(s)$ such that\break $\E
(\frac1n A_{cn/2}) \to\frac12 \int_0^c g(s)\,ds$.
\end{lemma}
\begin{pf}
Since $g_n(s) \le1$ almost surely, it suffices to show (by Fubini's
theorem and Lebesgue's dominated convergence theorem) that $\E
(g_n(sn/2)) \to g(s)$ for all fixed $s>0$. Let $\hat C_s$ be the cycle
of $\sigma_s$ containing the origin. By exchangeability, note that
\[
\E\bigl(g_n(sn/2)\bigr) = \P(v \in\hat C_{sn/2}),
\]
where $v$ is chosen uniformly among the $2L-1$ neighbors of 0. Fix such
a~neighbor~$v$. The idea for the proof of this lemma is that the cycle
structure of $v$ can be coupled with the cycle structure of the origin
in a random transposition process on the infinite line $\mathbb{Z}$,
rather than on the torus. More precisely, let $G_\infty$ be the graph
where the vertex set is $V_\infty= \mathbb{Z}$ and the edge set is
$E_\infty= \{(i,j) \in\mathbb{Z} \times\mathbb{Z}, |i-j| \le L\}$.
Consider the process \mbox{$(\sigma^\infty_t, t\ge0)$}, with values in the
permutation of $V_\infty$, obtained by transposing each edge \mbox{$(i,j)
\in E_\infty$} at rate $1/(2L)$. It is not obvious that this process is
well defined as there are an infinite number of edges. However, the
process may be constructed using a standard \textit{graphical
construction} [see, e.g., \citet{liggett}]. Briefly speaking, for every
(nonoriented) edge $e \in E_\infty$, consider an independent Poisson
process which rings at rate $1/(2L)$. Then the value $\sigma^\infty
_t(w)$ is defined for every $t\ge0$ and $w \in V_\infty$ by following
the trajectory between times~0 and $t$ of a particle which is initially
on $w$ and moves to a~neighbor~$j$ of its current position $i$ each
time the edge $e= (i,j)$ rings. It is easy to see (and will be shown
below) that almost surely there are empty patches (where no edge has
rung) surrounding the origin. Thus the trajectory cannot accumulate an
infinite number of jumps in a compact interval, and hence, is well defined.
Moreover, the cycle $\hat C^\infty_s$ of the origin in $\sigma
_s^\infty
$ contains only finitely many points almost surely for $s\ge0$, since
it must be contained in between two empty patches.

Let $c>0$. We claim that there is an event $\cG= \cG_n$ such that $\P
(\cG) \to1$ as $n \to\infty$ and such that on $\cG$, $\hat C_{sn/2}$
and $\hat C^\infty_{s}$ are identical. (Here we use the obvious
identification of $V = \mathbb{Z}/n\mathbb{Z}$ as a subset of
$\mathbb
{Z}$, as $V = \{-\lfloor n/2 \rfloor+1 , \ldots, -1, 0, 1, \ldots,
\lfloor n/2\rfloor\}$.) We choose $g(s) = \P( v \in\hat C^\infty_s)$.

The event $\cG$ we choose is
\[
\cG_n= \{\hat C_u \subset[ - \log n, \log n] \mbox{ for all } u \le
sn/2\}.
\]
The coupling between $\hat C_{sn/2}$ and $\hat C^\infty_s$ is obvious
on $\cG_n$ since we can use the same graphical construction for both
$\sigma_t$ and $\sigma^\infty_t$. It remains to show that $\P(\cG)
\to 1. $\vadjust{\goodbreak}
To do this it suffices that there is a strip of size at least $L$ in
$[-\log n,0]$ and in $[0,\log n]$ where each vertex in the strip has
never been involved in a transposition by time $sn/2$ (we say that such
a vertex has degree 0), what we called earlier an empty patch.
A given interval of size~$L$ contains exactly $L(2L-1) - {L \choose2}
= (3L^2 -L)/2 $ distinct edges, hence, the probability that it is an
empty patch is
\[
\exp\biggl(- \frac{s}{2L} \frac{3L^2 - L}2 \biggr) =:p(s) >0.
\]
If some patches of size $L$ share no edge in common, then the events
that they are empty are mutually independent. Since we can find at
least $ \alpha\log n$ distinct patches that do not share any edge in
$[0, \log n]$, for some $\alpha>0$ depending only on $L$, the
probability that there is no empty patch in $[0, \log n]$ is at most
$(1-p(s))^{\alpha\log n} \to0$. Hence, $\P(\cG_n) \to1$ and Lemma
\ref{LEA} is proved.
\end{pf}
%
%le19 #&#
\begin{lemma}\label{LvarA}
$\var(\frac1n A_{cn/2} ) \to0$ as $n \to\infty$.
\end{lemma}
\begin{pf}
Using (\ref{Echgvar}) and Cauchy--Schwarz's inequality,
\[
\var\biggl(\frac1n A_{cn/2}\biggr) \le\frac{c}4 \int_0^c
\var\bigl(g_n(sn/2)\bigr)\,ds.
\]
Since $g_n(s) \le1$, it suffices to show that $\var(g_n(sn/2)) \to0$
for all fixed $s>0$. Now, note that
\[
g_n(sn/2) = \frac1n\sum_{v \in V} f_v,
\]
where
\[
f_v = \frac1{2L-1} \sum_{\|w - v\| \le L} \mathbf{1}_{\{w \in\hat
C_v(sn/2)\}},
\]
and where $\hat C_v(s)$ denotes the cycle containing $v$ in $\sigma
_s$. Let
\[
\cA_v = \{\hat C_v(r) \subset[ v - \log n, v+\log n] \mbox{ for all
} r
\le sn/2\},
\]
where the addition and substraction is done modulo $n$. If $\| v - v'\|
> 2 \log n$, then on $\cA_v \cap\cA_{v'}$ the random variables $f_v$
and $f_{v'}$ may be taken to be independent. Reasoning as in Lemma \ref
{LconcX} shows that $\var(\frac1n \sum_v f_v) \to0$, since by Lemma
\ref{LEA} we know that $\P(\cA_v \cap\cA_{v'}) \to1$.
\end{pf}

Our final step is to show that $M_{cn/2}/n$ converges in probability to 0.
%
%le20 #&#
\begin{lemma}
\label{LconcM}
For all $\eps>0$,
\[
\P\Bigl( {\sup_{s \le cn/2}} |M_s| > \eps n\Bigr) \to0.
\]
\end{lemma}
\begin{pf}
By Markov's inequality,
\begin{eqnarray*}
\P\Bigl({\sup_{s\le t}} |M_s/n|>\eps\Bigr)&=&\bbp\Bigl({\sup_{s\le
t}}|M_s/n|^2>\eps^2\Bigr)\\
&\le&\frac{\E({\sup_{s\le t}}|M_s/n|^2)}{\eps^2}\\
&\le&\frac{4\E(M_t^2/n^2)}{\eps^2}
\end{eqnarray*}
by Doob's inequality. Now note that since $M_t$ is a martingale whose
jumps are only of size 1,
\[
M_t^2 - \int_0^t g_n(s)\,ds
\]
is again an $(\cF_t)_{t\ge0}$-martingale. [To see this, observe that
$F_{A_t^{-1}}$ is Poi($t$) and hence, $M^2_{A_t^{-1}}-t$ is a
martingale.] Thus $\E(M_t^2) \le t$ for all $t\ge0$ and when $t = cn/2$,
\[
\P\Bigl( {\sup_{s \le cn/2}} |M_s| > \eps n\Bigr) \le\frac{2c}{n
\eps^2} \to0
\]
as claimed.
\end{pf}
\begin{pf*}{Proof of Theorem \ref{Tmicro}}
It follows from Lemmas \ref{LEA} and \ref{LvarA} that $\frac1n
A_{cn/2} \to\frac12 \int_0^c g(s)\,ds$ in probability, where $g(s) =
\P(v \in\hat C^\infty_s)$ has been defined in Lem\-ma~\ref{LEA}. By
Lemma \ref{LconcM}, we deduce that
\[
\frac1n F_{cn/2} \to\frac1{2}\int_0^c g(s)\,ds
\]
in probability.
Since $\delta(t) = N_t - 2F_t$ for all $t \ge0$, it follows that
\[
\frac1n \delta(cn/2) \to v(c) = \frac{c}2 - \int_0^c g(s) \,ds.
\]
By Lemma \ref{Lboundedfrags}, we must have $v(c) < c/2$ for all $c>0$.
It thus suffices to show that $g$ is continuously differentiable on
$[0, \infty)$.
Assume that the process $(\sigma^\infty_t)$ is in some state such that
the (finite) cycle $C$ containing 0 also contains~$v$. Let~$f_1(C)$
denote the instantaneous rate at which $v$ becomes part of a~different
cycle; note that this rate depends indeed only on $C$ and not on the
rest of $\sigma^\infty_t$, and satisfies $f_1(C) \le|C|^2 / (2L)$.
Likewise, assume that the cycle containing $v$, $C'$, is distinct from
$C$. Let $f_2(C,C')$ be the instantaneous rate at which these cycles
merge. Then $f_2(C, C') \le|C| \times|C'| / (2L)$.

Note that $|\hat C_c^\infty|\ge k$ implies that there are $\lfloor
k/L\rfloor$ consecutive intervals of size $L$ around 0 all containing
at least one edge in the associated percolation process. By considering
every other interval, this implies that we can find $\lfloor
k/(2L)\rfloor$ disjoint intervals of size $L$, all of which contain at
least one edge. Such events are independent, and hence, if $p_\infty
(c)>0$ is the probability that at time $c$ an interval of size $L$ is
an empty patch, we find (summing over at most $k$ possible locations
for the leftmost point of this sequence of consecutive intervals),
%$$g(c+h) - g(c) = \E[\indic{v\in\hat C^\infty_c} hf_1(\hat C^
%(\hat%C_c^\infty, \hat C_c^\infty(v))+o(h)],$$
%where $\hat C_c^\infty(v)$ is the cycle of $v$ in $\sigma_c^\infty$.
%If we subtract $g(c)$ and divide by $h$ on the right-hand side above,
%then the limit as $h\to0$ exists from the following observations. We
%have
%
\[
\P(|\hat C_c^\infty|\ge k) \le k\bigl(1-p_\infty(c)\bigr)^{\lfloor k/2L
\rfloor},
\]
so that $|C_c^\infty|$ has exponential tails. It follows directly that
$\E(| \hat C_c^\infty|^2) < \infty$, and if $C_c^\infty(v)$ denotes the
cycle containing $v$ at time $c$, $\E(|\hat C_c^\infty||\hat
C_c^\infty
(v)|) < \infty$ by Cauchy--Schwarz's inequality. A routine argument
thus shows that
\[
g'(c) = \E\bigl[\mathbf{1}_{\{v\notin\hat C_c^\infty\}} f_2
(\hat
C_c^\infty
, \hat C_c^\infty(v))\bigr]-\E\bigl[\mathbf{1}_{\{v\in
\hat C^\infty_c\}}
f_1(\hat C^\infty_c) \bigr].
\]
By the same arguments, we see that $g'(c) $ is continuous, which in
turn shows that $v$ is continuously twice differentiable. The proof of
Theorem \ref{Tmicro} is complete.
\end{pf*}

%suskaldyti doi

% imsref loaded by lrinkeviciute, 2012-02-28 07:48:03

\printaddresses

\end{document}